\magnification=\magstep1
\input amstex
\documentstyle{amsppt}

\define\defeq{\overset{\text{def}}\to=}
\define\ab{\operatorname{ab}}
\define\pr{\operatorname{pr}}
\define\Gal{\operatorname{Gal}}
\define\Ker{\operatorname{Ker}}
\define\Ind{\operatorname{Ind}}
\define\et{\operatorname{et}}

\define\id{\operatorname{id}}
\define\cn{\operatorname{cn}}
\def \isom {\overset \sim \to \rightarrow}
\define\Sect{\operatorname{Sect}}
\define\Spec{\operatorname{Spec}}
\def\index{\operatorname{index}}
\def\cl{\operatorname{cl}}

\def\Spf{\operatorname{Spf}}
\def\rig{\operatorname{rig}}
\def\Fr{\operatorname{Fr}}
\def\Pic{\operatorname{Pic}}

\def\Sp{\operatorname{Sp}}

\def\geo{\operatorname{geo}}
\def\sol{\operatorname{sol}}
\def\nor{\operatorname{nor}}

\def\reg{\operatorname{reg}}
\def\red{\operatorname{red}}

\NoRunningHeads
\NoBlackBoxes
\topmatter

\title
Local sections of arithmetic fundamental groups of $p$-adic curves
\endtitle

\author
Mohamed Sa\"\i di
\endauthor
\abstract We investigate {\bf sections} of the arithmetic fundamental group $\pi_1(X)$ where $X$ is either a {\bf smooth affinoid $p$-adic curve}, 
or a {\bf formal germ of a $p$-adic curve}, and prove that they can be lifted (unconditionally) to sections of cuspidally abelian Galois groups.  As a consequence, if $X$ admits 
a compactification $Y$, and the exact sequence of $\pi_1(X)$ {\bf splits}, then $\text {index} (Y)=1$. 
We also exhibit a necessary and sufficient condition for a section of $\pi_1(X)$ to arise from a {\bf rational point} of $Y$.  
One of the key ingredients in our investigation is the fact, we prove in this paper in case $X$ is affinoid,
that the Picard group of $X$ is {\bf finite}.
\endabstract

\toc

\subhead
\S0. Introduction/Main results
\endsubhead

\subhead
\S1. Geometrically abelian arithmetic fundamental groups 
\endsubhead

\subhead
\S2. Cuspidally abelian arithmetic fundamental groups 
\endsubhead

\subhead
\S3. Geometric sections of arithmetic fundamental groups 
\endsubhead

\subhead
\S4. Picard groups of affinoid $p$-adic curves
\endsubhead

\subhead
\S 5. Compactification of formal germs of $p$-adic curves
\endsubhead

\endtoc

\endtopmatter

\document

\subhead
\S 0. Introduction/Main results
\endsubhead
This paper is motivated by the $p$-adic analog of the anabelian Grothendieck section conjecture.

Let $p\ge 2$ be a prime number, $k/\Bbb Q_p$ a finite extension, and $Y$ a proper, smooth, and geometrically connected hyperbolic $k$-curve. 
The arithmetic fundamental group $\pi_1(Y)$ of $Y$ projects onto the Galois group $G_k\defeq \Gal (\bar k/k)$ of $k$. A $k$-rational point
$x:\Spec k\to Y$ gives rise, by functoriality of fundamental groups, to a section $s_x:G_k\to \pi_1(Y)$ of the projection $\pi_1(Y)\twoheadrightarrow G_k$. 
We shall refer to such a section $s_x$ as {\bf geometric}.

\smallskip
\noindent
{\bf Question A.} {\it Is every section of the projection $\pi_1(Y)\twoheadrightarrow G_k$ geometric?}

\smallskip
In [Sa\"\i di2], Theorem 2 in the Introduction, we established two necessary and sufficient conditions for a group-theoretic section of the projection
$\pi_1(Y)\twoheadrightarrow G_k$ to be geometric. 
In [Hoshi] Hoshi  constructed a group-theoretic section $G_k\to \pi_1(Y)^{(p)}$ of the projection
$\pi_1(Y)^{(p)}\twoheadrightarrow G_k$ for a specific example $Y$, where $\pi_1(Y)^{(p)}$ is the geometrically pro-$p$ quotient of $\pi_1(Y)$,
which is {\bf not geometric} (i.e., doesn't arise from a scheme morphism $x:\Spec k\to Y$).
The author is not aware of any example of a $Y$ as above and a group-theoretic section of the projection $\pi_1(Y)\twoheadrightarrow G_k$ 
which is not geometric. 

Let $X$ be either a geometrically connected {\bf affinoid} subspace of $Y^{\rig}$, the rigid analytic curve associated to $Y$,
or a {\bf formal germ} of $Y$ meaning $X=\Spec (\hat {\Cal O}_{\Cal Y,y}\otimes _{\Cal O_k}k)$ is geometrically connected, where $\hat {\Cal O}_{\Cal Y,y}$ is the completion of the local ring
$\Cal O_{\Cal Y,y}$ of a model $\Cal Y$ of $Y$ over the ring of valuation $\Cal O_k$ of $k$ at a closed point $y\in \Cal Y^{\cl}$ (cf. Notations). 
Let $\pi_1(X)$ be the \'etale fundamental group of $X$ which sits in the exact sequence (cf. Notations)
$$1\to \pi_1(X)^{\geo}\to \pi_1(X)\to G_k\defeq \Gal (\bar k/k)\to 1.$$ 

A section $s:G_k\to \pi_1(X)$ of the projection $\pi_1(X)\twoheadrightarrow G_k$ induces a section $s_Y:G_k\to \pi_1(Y)$ of the projection $\pi_1(Y)\twoheadrightarrow G_k$
[cf. Notations, diagram (0.1)] which we shall refer to as a {\bf local section} of the projection $\pi_1(Y)\twoheadrightarrow G_k$. 
A geometric section is necessarily a local section as one easily verifies.
This prompts the following question, which motivates our study in this paper  of local sections of arithmetic fundamental groups of $p$-adic curves.

\smallskip
\noindent
{\bf Question B.} {\it Is every local section of the projection $\pi_1(Y)\twoheadrightarrow G_k$ geometric?}

\smallskip
Motivated by Questions A and B we investigate sections of arithmetic fundamental groups of affinoid $k$-curves and formal $p$-adic germs of curves.

Let $X$ be either a smooth and geometrically connected {\bf $k$-affinoid} curve or a {\bf formal $p$-adic germ} (cf. Notations for precise definitions).
Let $\pi_1(X)^{\geo,\ab}$ be the maximal abelian quotient of $\pi_1(X)^{\geo}$ and $\pi_1(X)^{(\ab)}$ the {\bf geometrically abelian} quotient of $\pi_1(X)$ 
which sits in the exact sequence
$$1\to \pi_1(X)^{\geo,\ab}\to \pi_1(X)^{(\ab)}\to G_k\to 1.$$ 

Similarly let $G_X\defeq \Gal (\overline L/L)$ be the absolute Galois group of the function field $L$ of $X$ (see Notations for the definition of $L$) which sits in the exact sequence (cf. $\S1$)
$$1\to G_X^{\geo}\to G_X\to G_k\to 1.$$
Let $G_X^{\geo,\ab}$
be the maximal abelian quotient of $G_X^{\geo}$ and $G_X^{(\ab)}$ the {\bf geometrically abelian} quotient of $G_X$ which sits in the exact sequence
$$1\to G_X^{\geo,\ab}\to G_X^{(\ab)}\to G_k\to 1.$$ 
We have an exact sequence 
$$1\to \widetilde {\Cal H}_X\to G_X^{(\ab)}\to \pi_1(X)^{(\ab)}\to 1$$ 
where $\widetilde {\Cal H}_X \defeq \Ker [G_X^{(\ab)}\twoheadrightarrow \pi_1(X)^{(\ab)}]$. 
In $\S1$ we investigate the structure of the $G_k$-module $\widetilde {\Cal H}_X$.
We prove in Proposition 1.4 that 
$\widetilde {\Cal H}_X$ is (canonically) isomorphic to 
$\prod _{x\in X^{\cl}} \Ind _{k(x)}^k \hat \Bbb Z(1)$ 
where the product is over all closed points of $X$
and $k(x)$ is the residue field at $x$ (cf. loc. cit.).

The Galois group $G_X$ sits in an exact sequence
$$1\to {\Cal H}_X\to G_X\to \pi_1(X)\to 1$$ 
where $\Cal H_X\defeq \Ker [G_X\twoheadrightarrow  \pi_1(X)]$.
Let $\Cal H_X^{\ab}$
be the maximal abelian quotient of $\Cal H_X$ and $G_X^{(\text c-\ab)}$ the {\bf geometrically cuspidally abelian} quotient of $G_X$ which sits in the exact sequence (cf. loc. cit.)
$$1\to \Cal H_X^{\ab}\to G_X^{(\text c-\ab)}\to \pi_1(X)\to 1.$$ 

In $\S2$ we investigate, in the framework of the theory of {\bf cuspidalisation} of sections of arithmetic fundamental groups (cf. [Sa\"\i di1] and [Sa\"\i di2]), 
sections $s:G_k\to \pi_1(X)$ of the projection $\pi_1(X)\twoheadrightarrow G_k$.
Let $Y$ be a {\bf $k$-compactification} of $X$ and $s_Y:G_k\to \pi_1(Y)$
the induced {\bf local section} of the projection $\pi_1(Y)\twoheadrightarrow G_k$ [cf. Notations for precise definitions and the diagram (0.1) therein].
One of our main results is the following [cf. Theorem 2.4 and Theorem 3.1(ii)].

\proclaim{Theorem A} {\bf (Lifting of sections to cuspidally abelian Galois groups)}. 
Let $s:G_k\to \pi_1(X)$ be a {\bf section} of the projection $\pi_1(X)\twoheadrightarrow G_k$. The followings hold.

(i)\ There {\bf exists} a section $s^{\text c-\ab}:G_k\to G_X^{(\text c-\ab)}$ of the projection $G_X^{(\text c-\ab)}\twoheadrightarrow G_k$ 
which {\bf lifts} the section $s$, i.e., which inserts in the following commutative diagram
$$
\CD
G_k @>s^{\text c-\ab} >> G_X^{(\text c-\ab)}\\
@| @VVV\\
G_k @>s>> \pi_1(X)\\
\endCD
$$
where the right vertical map is the natural projection $G_X^{(\text c-\ab)}\twoheadrightarrow \pi_1(X)$.
In particular, the set of sections of the 
projection $G_X^{(\text c-\ab)}\twoheadrightarrow G_k$ which lift the section $s$ is non-empty, and is (up to conjugation by elements of $\Cal H_X^{\ab}$) 
a torsor under $H^1(G_k,\Cal H_X^{\ab})$.

(ii)\ Assume $Y$ is hyperbolic. Then the section $s_Y:G_k\to \pi_1(Y)$ induced by $s$ is {\bf uniformly orthogonal to $\Pic$} 
in the sense of [Sa\"\i di1], Definition 1.4.1.
\endproclaim

The section $s$ is uniformly orthogonal to $\Pic$  (as in (ii) above) means that the retraction map $s^{*}:H^2(\pi_1(Y),\hat \Bbb Z(1))\isom H^2_{\et}(Y,\hat \Bbb Z(1))\to H^2(G_k,\hat \Bbb Z(1))$, which is induced by the section $s$, annihilates the Picard part of $H^2_{\et}(Y,\hat \Bbb Z(1))$, and similarly for every neighborhood $Y'\to Y$ of the section $s$.

Theorem A(ii) implies that local sections of arithmetic fundamental groups of hyperbolic $p$-adic curves satisfy condition (i) in 
[Sa\"\i di2], Theorem 2 in the Introduction.
In this sense local sections are close to being geometric. Establishing Theorem A(ii) was one of the main motivations for the author to investigate local sections of arithmetic fundamental groups of $p$-adic curves. Apart from local sections, and geometric sections, the author is not aware (for the time being) of any examples of group-theoretic sections of arithmetic fundamental groups of hyperbolic $p$-adic curves which are orthogonal to $\Pic$. 

As a consequence of Theorem A, and an observation of Esnault and Wittenberg on geometrically abelian sections of $p$-adic curves,
we deduce the following (cf. Theorem 2.5).

\proclaim {Theorem B} Assume that $X$ admits a {\bf $k$-compatctification} $Y$ (cf. Notations). If
the projection $\pi_1(X)\twoheadrightarrow G_k$ {\bf splits} then $\index (Y)=1$. 
\endproclaim

Theorem B asserts that the existence of local sections  of arithmetic fundamental groups of $p$-adic curves implies the existence of 
degree $1$ rational divisors. The link between sections of geometrically abelian Galois groups and the existence of 
degree $1$ rational divisors has been investigated in [Esnault-Wittenberg].

In $\S3$ we assume that $X$ admits a $k$-compactification $Y$ (cf. Notations). Let $\Pi_Y [X]$ be the \'etale fundamental group which 
classifies finite covers $Y'\to Y$ which only ramify at points of $Y$ {\bf not in} $X$ (cf. 3.3, as well as Notations for the meaning of ``not in $X$").
A section $s:G_k\to \pi_1(X)$ of the projection $\pi_1(X)\twoheadrightarrow G_k$ induces naturally a section $s^{\dagger}:G_k\to \Pi_Y[X]$ of the projection $\Pi_Y[ X]\twoheadrightarrow G_k$. We say that the section $s$ is {\bf geometric} (relative to $Y$) if the image $s^{\dagger}(G_k)$ is contained in a decomposition group 
$D_x\subset \Pi_Y[X]$ associated to a rational point $x\in Y(k)$ (cf. Definition 3.3.2). 
Further, we say that $s$ is {\bf admissible} (relative to $Y$) (cf. Definition 3.5.1) if for every open subgroup
$H\subset \Pi_Y[X]$ with $s^{\dagger}(G_k)\subset H$, corresponding to (a possibly ramified) cover $Y'\to Y$, 
the following holds. Let $G_{Y'}^{(1/p^2-\sol)}$ be the {\bf geometrically cuspidally $1/p^2$-solvable Galois group} of $Y'$: i.e., the maximal quotient $G_{Y'}\twoheadrightarrow H\twoheadrightarrow\pi_1(Y')$ of the absolute Galois group $G_{Y'}$ of $Y'$ such that $\Ker [H\twoheadrightarrow \pi_1(Y')]$ is abelian annihilated by $p^2$ (cf. [Sa\"\i di2], 3.1).
There exists a section $\tilde s_{Y'}:G_k\to G_{Y'}^{(1/p^2-\sol)}$ of the projection $G_{Y'}^{(1/p^2-\sol)}\twoheadrightarrow G_k$
[such a section exists unconditionally (see discussion in 3.5)] satisfying  the following property:

\smallskip
{\it For each open subgroup $F\subset G_{Y'}^{(1/p^2-\sol)}$ with $\tilde s_{Y'}(G_k)\subset F$, corresponding to a (possibly ramified) cover $Y''\to Y'$ with $Y''$ geometrically connected, the class of $\Pic^1_{Y''}$ in $H^1(G_k,\Pic^0_{Y''})$ is divisible by $p$}.

\smallskip
Our main result in $\S3$ is the following (cf. Theorem 3.5.2).

\proclaim {Theorem C} The section $s:G_k\to \pi_1(X)$ is {\bf geometric} (relative to $Y$)
{\bf if and only if} $s$ is {\bf admissible} (relative to $Y$).
\endproclaim

One of the key ingredients used in the proofs of the above results is the fact that $\Pic (X)$ is {\bf finite}.
In the case where $X$ is a {\bf formal $p$-adic germ} this is established in [Sa\"\i di2], Proposition 5.4, as a consequence of a result of Shuji Saito (cf. loc. cit.).
In case $X$ is {\bf affinoid} this is proven in $\S4$ (cf. Proposition 4.1) and may be of interest independently of 
the topics discussed in this paper. More precisely, we prove the following.

\proclaim{Theorem D (Picard groups of affinoid $p$-adic curves)}  Let $k$ be a {\bf $p$-adic local field} (i.e., $k/\Bbb Q_p$ is a finite extension) and $X=\Sp(A)$ a 
{\bf smooth} and geometrically connected  {\bf $k$-affinoid curve}. Then the Picard group $\Pic(X)$ is {\bf finite}.
\endproclaim

Finally in $\S5$ we prove (cf. Proposition 5.1) a compactification result for two dimensional complete local $p$-adic rings which is used in the proofs of
Proposition 1.2 and Proposition 2.2.

The results in $\S4$ and $\S5$ are used in this paper in $\S2$ and $\S3$, none of the results in  $\S2$ and $\S3$ is used in $\S4$ or $\S5$.

In this paper we worked with full arithmetic fundamental groups. Instead one could consider a similar setting and work with geometrically pro-$p$ arithmetic fundamental groups and Galois groups as in [Sa\"\i di2] (where one considers geometrically pro-$\Sigma$ arithmetic fundamental groups and Galois groups, 
$\Sigma$ being a set of primes containing $p$). In this geometrically pro-$p$ (pro-$\Sigma$) setting one can prove analogs of Theorems A and C.

\bigskip
\subhead
Acknowledgment
\endsubhead
I would like to thank Akio Tamagawa for several discussions we had on the topic of this paper. I would like to thank the referee for his/her careful reading of the paper and comments.

\bigskip
\subhead
Notations
\endsubhead
The following notations will be used throughout this paper (unless we specify otherwise).

\smallskip
$\bullet$  $p\ge 2$ is a prime number, $k$ is a {\bf $p$-adic local field} (i.e., $k/\Bbb Q_p$ is a finite extension) with ring of valuation $\Cal O_k$, uniformizer
$\pi$, and residue field $F$. Thus $F$ is a finite field of characteristic $p$.

\smallskip
$\bullet$ A proper, smooth, and geometrically connected $k$-curve $Y$ is {\bf hyperbolic} if $\text {genus} (Y)\ge 2$.

\smallskip
$\bullet$  For a profinite group $H$ we denote by $H^{\ab}$ the maximal {\bf abelian} quotient of $H$.

\smallskip
$\bullet$ Let 
$$1\to H'\to H @>{\pr}>> G \to 1$$ 
be an exact sequence of profinite groups. We will refer to
a continuous homomorphism $s:G\to H$ such that $\pr \circ s=\id_{G}$ as a (group-theoretic)
{\bf section} of the above sequence, or simply a section of the projection $\pr : H\twoheadrightarrow G$.

\smallskip
$\bullet$  All scheme cohomology groups considered in this paper are \'etale cohomology groups.

\subhead
Affinoid $p$-adic curves
\endsubhead

\smallskip
$\bullet$ $X=\Sp A$ is a {\bf smooth} and geometrically connected  {\bf affinoid $k$-curve}. On occasions we will write, if there is no risk of confusion, 
$X=\Spec A$ for the corresponding affine $k$-scheme.

\smallskip
$\bullet$ One can embed $X$ into a proper, smooth, and geometrically connected rigid analytic curve $Y^{\rig}: X\hookrightarrow Y^{\rig}$ 
so that $X$ is an open affinoid subspace of $Y^{\rig}$
(cf. [Fresnel-Matignon], 2.6, Corollaire 2).
Write $Y$ for the algebraization of $Y^{\rig}$ via the rigid GAGA functor 
which is a proper, smooth, and geometrically connected algebraic $k$-curve. We will refer to 
$X$ as a {\bf $p$-adic affinoid curve} (or simply an affinoid) and $Y$ a {\bf $k$-compactification} of $X$.

\subhead
Formal $p$-adic germs
\endsubhead

\smallskip
$\bullet$ $A$ is a {\bf normal two dimensional complete local ring} containing $\Cal O_k$ with maximal ideal $\frak m_A$ containing $\pi$ and residue field
$F=A/\frak m_A$. Write $A_k\defeq A\otimes _{\Cal O_k}k=A[\frac {1}{\pi}]$ and $X\defeq \Spec A_k$. 
We assume $X$ is geometrically connected and refer to $X$ as a {\bf formal $p$-adic germ}.

\smallskip
$\bullet$ A {\bf ($k$-)compactification} of $\Spec A$ is a proper and flat relative $\Cal O_k$-curve $\Cal Y\to \Spec \Cal O_k$ with $\Cal Y$ normal, 
$Y\defeq \Cal Y\times _{\Spec \Cal O_k}\Spec k$ geometrically connected, $y\in \Cal Y^{\cl}$ is a closed point, $\Cal O_{\Cal Y,y}$ is the local ring of $\Cal Y$ at $y$, 
$\hat {\Cal O}_{\Cal Y,y}$ its completion, with an isomorphism $\hat {\Cal O}_{\Cal Y,y}\isom A$. We have a natural scheme morphism $X\to Y$. 
We shall refer to $Y$ as a {\bf $k$-compactification} of $X$. In $\S5$ we prove the existence  of such a compactification 
$X\to Y$ after possibly a finite extension of $k$ (cf. loc. cit., Proposition 5.1).

\smallskip
In what follows $X$ is either {\bf an affinoid $p$-adic curve} or a {\bf formal $p$-adic germ}.

\smallskip
$\bullet$\ We say that $X$ is {\bf hyperbolic} if there exists a finite extension $k'/k$ such that
$X_{k'}\defeq \Spec(A\otimes_k k')$ [resp. $X_{k'}\defeq \Sp (A\otimes_k k')$ if $X$ is affinoid]
possesses a $k'$-compactification $Y$ with $Y$ hyperbolic. 
There exists a finite extension $k'/k$ and a finite geometric \'etale cover $X'\to X_{k'}$ with $X'$ geometrically connected and hyperbolic.
This is Proposition 5.3 in case $X$ is a {\bf formal $p$-adic germ} and follows from [Sa\"\i di3], Theorem A, in case $X$ is {\bf affinoid}.

\smallskip
$\bullet$\ $\eta$ is a fixed choice of a geometric point of $X$ with values in its generic point. Thus $\eta$ determines algebraic closures $\bar k$, 
$\overline L$, of $k$, and $L\defeq \Fr(A)$; respectively. We have an exact sequence of fundamental groups
$$1@>>>\pi_1(X,\eta)^{\geo}@>>> \pi_1(X,\eta)@>>> G_k @>>> 1$$
where $\pi_1(X,\eta)$ is the \'etale fundamental group of $X$ 
with geometric point $\eta$ (cf. [Sa\"\i di3], 2.1,  for more details on the definition of 
$\pi_1(X,\eta)$ in case $X$ is an affinoid), 
$\pi_1(X,\eta)^{\geo}\defeq \Ker[\pi_1(X,\eta)\twoheadrightarrow  G_k]$, 
and $G_k\defeq \Gal (\bar k/k)$ is the absolute Galois group of $k$. 

\smallskip
In what follows $Y$ is a {\bf $k$-compactification} of $X$.

\smallskip
$\bullet$\  We have a commutative diagram of exact sequences 
of arithmetic fundamental groups 
$$
\CD
1@>>>\pi_1(X,\eta)^{\geo}@>>> \pi_1(X,\eta)@>>> G_k @>>> 1 \\
@.@VVV  @VVV @|\\
1@>>> \pi_1(Y_{\bar k},\bar \eta)@>>> \pi_1(Y,\eta)@>>> G_k@>>> 1\\
\endCD
\tag 0.1
$$
where $\pi_1(Y,\eta)$ [resp. $\pi_1(Y_{\bar k},\bar \eta)$]
is the \'etale fundamental group of $Y$ (resp. $Y_{\bar k}\defeq Y\times _{\Spec k}{\Spec \bar k}$)
with geometric point $\eta$ (resp. $\bar \eta$ which is induced by $\eta$).
In case $X$ is an {\bf affinoid} (resp. a {\bf formal $p$-adic germ}) the middle vertical map is induced by the rigid analytic morphism $X\to Y^{\rig}$ and the rigid GAGA functor
(resp. the scheme morphism $X\to Y$).

\smallskip
$\bullet$ \ We write $X^{\cl}$ (resp. $Y^{\cl}$) for the set of closed points of $X$ (resp. $Y$). 
For a closed point $x$ of $X$ (resp. $Y$) we write $k(x)$ for the residue field at $x$. Thus $k(x)$ is a finite extension of $k$.

\smallskip
$\bullet$ \ We say that $x\in Y^{\cl}$ is {\bf not} in X if $x$ is not in the image of the scheme morphism $X\to Y$ if $X$ is a {\bf formal $p$-adic germ} 
or $x\notin X^{\cl}$ in case $X$ is {\bf affinoid}. In case $X=\Spec (\Cal O_{\Cal Y,y}\otimes _{\Cal O_k}k)$ is a {\bf formal $p$-adic germ} the set of closed points of $Y$ {\bf not} in 
$X$ is in one-to-one correspondence with the set of closed points of $Y$ which {\bf do not} specialise in $y$ [cf. [Liu], $\S10$, Proposition 1.40(a)].

\bigskip
{\bf Throughout sections $\S1$, $\S2$, and $\S3$,  $X$ will denote either an affinoid $p$-adic curve or a formal $p$-adic germ}. {\bf In $\S3$ we will assume $X$ admits a $k$-compactification $Y$ which is hyperbolic and fix a choice of such a compactification throughout}.
\bigskip
\subhead
\S 1. Geometrically abelian arithmetic fundamental groups 
\endsubhead
In this section we investigate the structure of various geometrically abelian arithmetic fundamental groups and absolute Galois group associated to $X$.
Let 
$$\pi_1(X,\eta)^{(\ab)}\defeq \pi_1(X,\eta)/\Ker [\pi_1(X,\eta)^{\geo}\twoheadrightarrow \pi_1(X,\eta)^{\geo,\ab}]$$
be the {\bf geometrically abelian} fundamental group of $X$ [here $\pi_1(X,\eta)^{\geo,\ab}$ denotes the maximal abelian quotient of $\pi_1(X,\eta)^{\geo}$].

\proclaim {Proposition 1.1} We use the above notations. The followings hold.

(i) Assume $X$ is an {\bf affinoid}. For each prime number $\ell$ the pro-$\ell$-Sylow subgroup of $\pi_1(X,\eta)^{\geo,\ab}$ is pro-$\ell$ abelian {\bf free}, of infinite rank if $\ell=p$, and finite (computable) rank otherwise (see [Sa\"\i di3], Theorem A, for the precise value of this rank in case $\ell \neq p$). 

(ii) Assume $X$ is a {\bf formal $p$-adic germ}. For each prime number $\ell\neq p$ the pro-$\ell$-Sylow subgroup of $\pi_1(X,\eta)^{\geo,\ab}$ is pro-$\ell$ abelian {\bf free} of finite computable rank (see [Sa\"\i di4], Theorem A, for the precise value of this rank).
\endproclaim

\demo{Proof} Assertion (i) follows from [Sa\"\i di3], Theorem A. [Note that the assumption in loc. cit. that $X$ is the complement in a proper rigid analytic $k$-curve of the disjoint union of finitely many $k$-rational open discs is satisfied after a finite extension of $k$ (cf. [Fresnel-Matignon], 2.6, Th\'eor\`eme 6 and Corollaire 1)].
Assertion (ii) follows from [Sa\"\i di4], Theorem A.
\qed
\enddemo

Let $S\defeq \{x_1,\ldots,x_n\}\subset X^{\cl}$ be a finite set of closed points and write $U\defeq X\setminus S$ viewed as an open sub-scheme of $X$ (resp. $X=\Spec A$
in case $X$ is an affinoid). Let $\pi_1(U,\eta)$ be the \'etale fundamental group of
$U$ with geometric point $\eta$ (cf. [Sa\"\i di3], 2.1, for the definition of $\pi_1(U,\eta)$ in case $X$ is affinoid) which sits in the exact sequence
$$1@>>> \pi_1(U,\eta)^{\geo}@>>> \pi_1(U,\eta)@>>> G_k @>>> 1$$
where $\pi_1(U,\eta)^{\geo}\defeq \Ker[\pi_1(U,\eta) \twoheadrightarrow G_k]$ (cf. loc. cit. in case $X$ is affinoid).
Let
$$\pi_1(U,\eta)^{(\ab)}\defeq \pi_1(U,\eta)/\Ker [\pi_1(U,\eta)^{\geo}\twoheadrightarrow \pi_1(U,\eta)^{\geo,\ab}]$$
be the {\bf geometrically abelian} fundamental group of $U$ [here $\pi_1(U,\eta)^{\geo,\ab}$ is the maximal abelian quotient of $\pi_1(U,\eta)^{\geo}$]. 
We have an exact sequence
$$1\to \widetilde \Delta_U\to \pi_1(U,\eta)^{(\ab)} \to \pi_1(X,\eta)^{(\ab)} \to 1\tag 1.1$$
where $\widetilde \Delta_U\defeq \Ker [\pi_1(U,\eta)^{(\ab)} \twoheadrightarrow \pi_1(X,\eta)^{(\ab)}]=\Ker [\pi_1(U,\eta)^{\geo,\ab} \twoheadrightarrow \pi_1(X,\eta)^{\geo,\ab}]$
and the (surjective) map $\pi_1(U,\eta)^{(\ab)} \twoheadrightarrow \pi_1(X,\eta)^{(\ab)}$ 
is induced by the natural projection $\pi_1(U,\eta)\twoheadrightarrow \pi_1(X,\eta)$. Note that $\widetilde \Delta_U$ has a natural structure of $G_k$-module.

\proclaim {Proposition 1.2} We use the above notations. There exists a natural isomorphism
$$\prod_{i=1}^n \Ind _{k(x_i)}^k \hat \Bbb Z(1)\isom \widetilde \Delta_U$$ 
of $G_k$-modules where the $(1)$ is a Tate twist.
\endproclaim

\demo{Proof} We have a natural surjective homomorphism $\prod_{i=1}^n \Ind _{k(x_i)}^k \hat \Bbb Z(1)\twoheadrightarrow \widetilde \Delta_U$ of $G_k$-modules
 mapping $\Ind_{k(x_i')}^k \hat \Bbb Z(1)$ onto the inertia subgroup [of $\pi_1(U,\eta)^{(\ab)}$]
 at $x_i$, as follows from the structure of inertia groups of Galois extensions of henselian 
 discrete valuation rings of residue characteristic zero. We show this map is an isomorphism. 
To this end we can, without loss of generality, assume that $X$ admits a $k$-compactification $Y$ (cf. Notations). Indeed, this holds for $X$ affinoid (cf. loc. cit.), 
and holds after possibly replacing $k$ by a finite field extension in case $X$ is a formal $p$-adic germ (cf. Proposition 5.1) which doesn't alter the structure of $\widetilde \Delta_U$.
We have a commutative diagram of exact sequences 

$$
\CD
1@>>>\pi_1(X,\eta)^{\geo,\ab}@>>> \pi_1(X,\eta)^{(\ab)}@>>> G_k @>>> 1 \\
@.@VVV  @VVV @|\\
1@>>> \pi_1(Y_{\bar k},\bar \eta)^{\ab}@>>> \pi_1(Y,\eta)^{(\ab)}@>>> G_k@>>> 1\\
\endCD
\tag 1.2
$$
where $\pi_1(Y,\eta)^{(\ab)}\defeq  \pi_1(Y,\eta)/\Ker[\pi_1(Y_{\bar k},\bar \eta)\twoheadrightarrow \pi_1(Y_{\bar k},\bar \eta)^{\ab}]$
and the middle vertical map is induced by the natural homomorphism $\pi_1(X,\eta)\to \pi_1(Y,\eta)$ [cf. Notations, diagram (0.1)].

Denote by $x_i'$ the image of $x_i$ in $Y$, $\forall 1\le i\le n$ [note that $k(x_i)=k(x_i')$]. Let $x_0'\in Y^{\cl} \setminus \{x_1',\ldots,x_n'\}$ be a closed point which is not in the image of $X$ (cf. Notations). 
Write
$S'\defeq \{x_0',x_1',\ldots,x_n'\}\subset Y^{\cl}$ and $V\defeq Y\setminus S'$ which is an affine $k$-curve.
Let $\pi_1(V,\eta)$ be the \'etale fundamental group of
$V$ with geometric point $\eta$ which sits in the exact sequence
$1@>>> \pi_1(V_{\bar k},\bar \eta)@>>> \pi_1(V,\eta)@>>> G_k @>>> 1,$
where $\pi_1(V_{\bar k},\bar \eta)$ is the \'etale fundamental group of $V_{\bar k}\defeq V\times _{k}\bar k$ with geometric point $\bar \eta$ which is induced by $\eta$.
Let
$\pi_1(V,\eta)^{(\ab)}\defeq \pi_1(V,\eta)/\Ker [\pi_1(V_{\bar k},\bar \eta)\twoheadrightarrow \pi_1(V_{\bar k},\bar \eta)^{\ab}]$
be the geometrically abelian fundamental group of $V$. We have a commutative diagram of exact sequences
$$
\CD
1@>>>\widetilde \Delta_U@>>> \pi_1(U,\eta)^{(\ab)}@>>>\pi_1(X,\eta)^{(\ab)}  @>>> 1 \\
@.@VVV  @VVV @VVV\\
1@>>> \widetilde \Delta_V@>>> \pi_1(V,\eta)^{(\ab)}@>>> \pi_1(Y,\eta)^{(\ab)} @>>> 1\\
\endCD
\tag 1.3
$$
where $\widetilde \Delta_V\defeq \Ker[\pi_1(V,\eta)^{(\ab)}\twoheadrightarrow \pi_1(Y,\eta)^{(\ab)}]$. 
The middle vertical map in diagram (1.3) is induced by the natural homomorphism $\pi_1(U,\eta)\to \pi_1(V,\eta)$, which is induced
by the scheme morphism $X\to Y$ in case $X$ is a formal $p$-adic germ, and
by the rigid analytic morphism $X\to Y^{\rig}$ and the rigid GAGA functor in case $X$ is affinoid (here we use the fact that $x_0'$ is not in the image of $X$). 
The right vertical map in diagram (1.3) is the middle vertical map in diagram (1.2).

One has an exact sequence of $G_k$-modules (as follows from the well-known structure of $\pi_1(V,\eta)^{(\ab)}$, see for example
the discussion in [Sa\"\i di5], $\S0$)
$$0\to \hat \Bbb Z(1)\to  \prod_{i=0}^n\Ind_{k(x_i')}^k \hat \Bbb Z(1)\to \widetilde \Delta_V\to 0.$$
Consider the composite homomorphism $\tau:\prod_{i=1}^n\Ind_{k(x_i')}^k \hat \Bbb Z(1)@>>> \widetilde \Delta_V$
of $G_k$-modules:
$$\prod_{i=1}^n\Ind_{k(x_i')}^k \hat \Bbb Z(1)\hookrightarrow \prod_{i=0}^n\Ind_{k(x_i')}^k \hat \Bbb Z(1)\twoheadrightarrow \widetilde \Delta_V$$
where the first map is the natural embedding: $(\beta_1,\ldots,\beta_n)\mapsto (0,\beta_1,\ldots,\beta_n)$ and the second map is as in the above exact sequence. 
Thus $\tau$ is injective (cf. above exact sequence). 
Consider the following commutative diagram
$$
\CD
\prod_{i=1}^n\Ind_{k(x_i')}^k \hat \Bbb Z(1) @>>> \widetilde \Delta_U\\
@VVV     @VVV\\
\prod_{i=0}^n\Ind_{k(x_i')}^k \hat \Bbb Z(1) @>>> \widetilde \Delta_V\\
\endCD
$$
where the right vertical map is the one in diagram (1.3). The left vertical and lower horizontal maps are as explained above, hence their composite is the map $\tau$.
The upper horizontal map is the natural projection $\prod_{i=1}^n\Ind_{k(x_i')}^k \hat \Bbb Z(1) \twoheadrightarrow \widetilde \Delta_U$ mentioned at the start of the proof. This map is an isomorphism since it is onto and it is injective as $\tau$ is.
\qed
\enddemo

\definition {Remark 1.3} With the notations in Proposition 1.2 and the proof therein assume that $x'_0\in Y(k)$ is a $k$-rational point.
In this case $\tau  (\prod_{i=1}^n\Ind_{k(x_i')}^k \hat \Bbb Z(1) )=\widetilde \Delta_V$, the 
map $\widetilde \Delta_U\to \widetilde \Delta_V$ is an isomorphism, and the right square in diagram (1.3) (cf. proof of Proposition 1.2) is cartesian. 
\enddefinition

Let $G_X\defeq \Gal (\overline L/L)$ [recall $L\defeq \Fr(A)$] which sits in the exact sequences 
$$1\to G_{X}^{\geo}\to G_X\to G_k\to 1$$
where $G_{X}^{\geo}\defeq \Gal (\overline L/L\bar k)$, and
$$1\to \Cal H_X\to G_X\to \pi_1(X,\eta)\to 1\tag 1.4$$
where $\Cal H_X\defeq \Ker [G_X\twoheadrightarrow \pi_1(X,\eta)]$. 
Let 
$$G_X^{(\ab)}\defeq G_X/\Ker (G_{X}^{\geo}\twoheadrightarrow G_{X}^{\geo,\ab})$$ 
which we shall refer to as the {\bf geometrically abelian} Galois group of $X$ (here $G_{X}^{\geo,\ab}$ is the maximal abelian quotient of $G_{X}^{\geo}$).
We have an exact sequence
$$1\to \widetilde {\Cal H}_X\to G_X^{(\ab)}\to \pi_1(X,\eta)^{(\ab)}\to 1\tag 1.5$$
where $\widetilde {\Cal H}_X\defeq \Ker [G_X^{(\ab)}\twoheadrightarrow \pi_1(X,\eta)^{(\ab)}]=\Ker [G_{X}^{\geo,\ab}\twoheadrightarrow \pi_1(X,\eta)^{\geo,\ab}]$.
Note that $\widetilde {\Cal H}_X$ has a natural structure of $G_k$-module.

\proclaim {Proposition 1.4} We use the above notations. There exists a natural isomorphism of $G_k$-modules
$$\prod _{x\in X^{\cl}} \Ind _{k(x)}^k \hat \Bbb Z(1)\isom \widetilde {\Cal H}_X$$
where the product is over all closed points $x\in X^{\cl}$.
\endproclaim

\demo{Proof} This follows from Proposition 1.2 and the fact that $\widetilde {\Cal H}_X\isom \varprojlim _{U}\widetilde \Delta _U$ where $U=X\setminus S$; $S$ runs over all finite
subsets of $X^{\cl}$, and $\widetilde \Delta_U$ is as in the proof of loc. cit.. (Note that  $G_X^{(\ab)}\isom \varprojlim _{U} \pi_1(U,\eta)^{(\ab)}$ where the limit runs over all $U$ as above.)
\qed
\enddemo

\subhead
\S 2. Cuspidally abelian arithmetic fundamental groups
\endsubhead
In this section we investigate the problem of {\bf cuspidalisation} of sections of the projection $\pi_1(X,\eta)\twoheadrightarrow G_k$.
This problem has been investigated in the case of proper and smooth hyperbolic $p$-adic curves in [Sa\"\i di1] and [Sa\"\i di2].
We use the notations in $\S0$ and $\S1$. 

Let $S\defeq \{x_1,\ldots,x_n\}\subset X^{\cl}$ be a finite set of closed points and $U\defeq X\setminus S$ (cf. $\S1$). 
Consider the exact sequence 
$$1\to\Delta _U \to \pi_1(U,\eta)^{\geo}\to \pi_1(X,\eta)^{\geo}\to 1$$
where $\Delta _U\defeq \Ker [ \pi_1(U,\eta)^{\geo}\twoheadrightarrow \pi_1(X,\eta)^{\geo}]$.
The maximal abelian quotient $\Delta_U^{\ab}$ of $\Delta_U$ is a $\pi_1(X,\eta)^{\geo}$-module.
Let $\Delta_U^{\cn}$ be the maximal quotient of $\Delta _{U}^{\ab}$ on which $\pi_1(X,\eta)^{\geo}$ acts trivially. Define 
$$\pi_1(U,\eta)^{\geo,\text c-\ab}\defeq \pi_1(U,\eta)^{\geo}/\Ker (\Delta_U\twoheadrightarrow \Delta_U^{\ab})$$
and 
$$\pi_1(U,\eta)^{\geo,\text c-\cn}\defeq \pi_1(U,\eta)^{\geo}/\Ker (\Delta _U   \twoheadrightarrow \Delta_U^{\cn}).$$

We shall refer to $\pi_1(U,\eta)^{\geo,\text c-\ab}$ (resp. $\pi_1(U,\eta)^{\geo,\text c-\cn}$)
as the {\bf cuspidally abelian} (resp. {\bf cuspidally central}) quotient of $\pi_1(U,\eta)^{\geo}$. 
Further, define
$$\pi_1(U,\eta)^{(\text c-\ab)}\defeq \pi_1(U,\eta)/\Ker (\Delta_U\twoheadrightarrow \Delta_U^{\ab})$$
and 
$$\pi_1(U,\eta)^{(\text c-\cn)}\defeq \pi_1(U,\eta)/\Ker  (\Delta_U   \twoheadrightarrow \Delta_U^{\cn}).$$
We shall refer to $\pi_1(U,\eta)^{(\text c-\ab)}$ (resp. $\pi_1(U,\eta)^{(\text c-\cn)}$)
as the {\bf (geometrically) cuspidally abelian} [resp. {\bf (geometrically) cuspidally central}] quotient of $\pi_1(U,\eta)$.
We have the following commutative diagram of exact sequences

$$
\CD
1@>>>\Delta_U@>>> \pi_1(U,\eta) @>>> \pi_1(X,\eta) @>>> 1 \\
@.@VVV  @VVV @|\\
1@>>> \Delta _U^{\ab}@>>> \pi_1(U,\eta)^{(\text c-\ab)} @>>> \pi_1(X,\eta) @>>> 1\\
@.@VVV  @VVV @|\\
1@>>> \Delta _U^{\cn}@>>> \pi_1(U,\eta)^{(\text c-\cn)} @>>> \pi_1(X,\eta) @>>> 1\\
@.@VVV  @VVV @VVV\\
1@>>> \widetilde \Delta_U@>>> \pi_1(U,\eta)^{(\ab)} @>>> \pi_1(X,\eta)^{(\ab)} @>>> 1
\endCD
\tag2.1
$$
where the middle vertical maps are surjective, and the
middle vertical map in the lower diagram is induced by the natural surjective map $\pi_1(U,\eta)^{\geo,\text c-\ab}\twoheadrightarrow  \pi_1(U,\eta)^{\geo,\ab}$. 
[Note that $\pi_1(X,\eta)^{\geo}$ acts trivially on the quotient $\widetilde \Delta_U$ of $\Delta_U^{\ab}$.]

\proclaim {Lemma 2.1} We use the above notations. The homomorphism $\Delta_U^{\cn}\to  \widetilde \Delta_U$ in diagram 
(2.1) is an {\bf isomorphism} of $G_k$-modules. In particular, the lower right square in diagram (2.1) is cartesian.
\endproclaim

\demo{Proof} The proof follows from Proposition 1.2 and the various definitions. More precisely, there exists a natural surjective homomorphism
$\prod_{i=1}^n\Ind_{k(x_i)}^k \hat \Bbb Z(1)\twoheadrightarrow \Delta_U^{\cn}$ 
[mapping $\Ind_{k(x_i')}^k \hat \Bbb Z(1)$ onto the inertia subgroup of $\pi_1(U,\eta)^{(\text c-\cn)}$
 at $x_i$, as follows from the structure of inertia groups of Galois extensions of henselian discrete valuation rings of residue characteristic zero]
which composed with the projection $\Delta _U^{\cn} \twoheadrightarrow
\widetilde \Delta_U$ is the isomorphism $\prod_{i=1}^n \Ind _{k(x_i)}^k \hat \Bbb Z(1)\isom \widetilde \Delta_U$ in Proposition 1.2 hence our assertion.
\qed
\enddemo

Let $s:G_k\to \pi_1(X,\eta)$ be a {\bf section} of the projection $\pi_1(X,\eta)\twoheadrightarrow G_k$. 

\proclaim{Proposition 2.2} {\bf (Lifting of sections to cuspidally central arithmetic fundamental groups)}.  
We use the above notations. There exists a section $s_U^{\text c-\cn}:G_k\to \pi_1(U,\eta)^{(\text c-\cn)}$ of the projection $\pi_1(U,\eta)^{(\text c-\cn)}\twoheadrightarrow G_k$ 
which {\bf lifts} the section $s$, i.e., which inserts in the following commutative diagram
$$
\CD
G_k @>s_U^{\text c-\cn} >> \pi_1(U,\eta)^{(\text c-\cn)}\\
@| @VVV\\
G_k @>s>> \pi_1(X,\eta)\\
\endCD
$$
where the right vertical map is the natural projection $\pi_1(U,\eta)^{(\text c-\cn)} \twoheadrightarrow \pi_1(X,\eta)$.
In particular, the set of sections of the 
projection $\pi_1(U,\eta)^{(\text c-\cn)}\twoheadrightarrow G_k$ which lift the section $s$ is non-empty, and is (up to conjugation by elements of $\Delta_U^{\cn}$) a torsor under $H^1(G_k,\Delta_U^{\cn})$.
\endproclaim

\demo{Proof} Consider the commutative diagram of exact sequences

$$
\CD
1@>>> \Delta_U^{\cn}@>>> E_U\defeq E_U[s] @>>> G_k @>>> 1\\
@. @| @VVV @VsVV\\
1@>>> \Delta_U^{\cn}@>>> \pi_1(U,\eta)^{(\text c-\cn)} @>>> \pi_1(X,\eta) @>>> 1
\endCD
$$
where the right square is cartesian. Thus the group extension $E_U$ is the pull-back of the group extension $\pi_1(U,\eta)^{(\text c-\cn)}$ by the section $s$. The set of 
(possible) splittings of the group extension $E_U$ is in one-to-one correspondence with the set of sections of the projection $\pi_1(U,\eta)^{(\text c-\cn)}\twoheadrightarrow G_k$ 
which lift the section $s$. We show the group extension $E_U$ splits.

To this end we can replace $k$ by a finite extension over which the points $\{x_i\}_{i=1}^n$ are rational, and we can also assume $n=1$ (see the argument at the start 
of the proof of Lemma 2.3.1 in [Sa\"\i di1]). 
Further we can replace $X$ by a neighbourhood $X'$ of the section $s$: i.e., an \'etale cover $X'\to X$ corresponding to an open subgroup $H=\pi_1(X',\eta)\subset \pi_1(X,\eta)$ containing the image $s(G_k)$ of $s$. Indeed, if $U'\defeq U\times _XX'$ there exists a commutative diagram of natural homomorphisms
$$
\CD
\pi_1(U',\eta)^{(\text c-\cn)}     @>>>  \pi_1(U,\eta)^{(\text c-\cn)}\\ 
@VVV @VVV\\
\pi_1(X',\eta)     @>>>  \pi_1(X,\eta)\\ 
\endCD
$$
where the upper horizontal map is induced by the natural map $\pi_1(U',\eta) \to  \pi_1(U,\eta)$ [note $\Delta_{U'}=\Delta_U$ and $\pi_1(X',\eta)^{\geo}$ acts trivially 
on $\Delta_U^{\cn}$], 
and the various maps in this diagram commute with the
projections onto $G_k$. 
The section $s$ induces a section $\tilde s:G_k\to \pi_1(X',\eta)$ of the projection $\pi_1(X',\eta)\twoheadrightarrow G_k$, 
and a lifting $\tilde s_{U'}^{\text c-\cn}:G_k\to \pi_1(U',\eta)^{(\text c-\cn)}$ of $\tilde s$ (as in the statement of Proposition 2.2) 
induces a lifting $s_U^{\text c-\cn}:G_k\to \pi_1(U,\eta)^{(\text c-\cn)}$ of $s$ as required (cf. above diagram).
Now it follows from [Sa\"\i di3], Theorem A, in case $X$ is an affinoid, and Proposition 5.3 in this paper (cf. $\S5$)
in case $X$ is a formal $p$-adic germ,
that there exists (after possibly a finite extension of $k$)
a neighbourhood $X'\to X$ of $s$ with $X'$ hyperbolic (cf. Notations).
We can thus assume, without loss of generality, that $X$ possesses a $k$-compactification $Y$ with $Y$ hyperbolic and 
the set $S\defeq \{x\}\subset X(k)$ consists of a single $k$-rational point, in which case $\Delta_U^{\cn}\isom \hat \Bbb Z(1)$ as a $\pi_1(X,\eta)$-module 
(cf. Lemma 2.1 and Proposition 1.2).



Consider the following maps (here $X=\Spec A$ in case $X$ is {\bf affinoid})
$$H^2(\pi_1(X,\eta),\hat \Bbb Z(1)) \hookrightarrow H^2(X,\hat \Bbb Z(1)) \leftarrow \Pic (X)$$
where the map $H^2(\pi_1(X,\eta),\hat \Bbb Z(1)) \hookrightarrow H^2(X,\hat \Bbb Z(1))$ arises from the Cartan-Leray spectral sequence
and is injective (cf. [Serre], proof of Proposition 1), and the map $\Pic(X)\rightarrow H^2(X,\hat \Bbb Z(1))$
is the cycle class map arising from the Kummer exact sequence in \'etale topology. Let $[\pi_1(U,\eta)^{(\text c-\cn)}]\in  H^2(\pi_1(X,\eta),\hat \Bbb Z(1))$
be the class of the group extension $\pi_1(U,\eta)^{(\text c-\cn)}$.
The image of $[\pi_1(U,\eta)^{(\text c-\cn)}]$ in $H^2(X,\hat \Bbb Z(1))$ 
coincides with the image of the line bundle $\Cal O(x)\in \Pic(X)$ via the Kummer map $\Pic (X) \rightarrow H^2(X,\hat \Bbb Z(1))$.
Indeed, this follows from the following commutative diagram
$$
\CD
H^2(\pi_1(X,\eta),\hat \Bbb Z(1)) @>>> H^2(X,\hat \Bbb Z(1)) @<<< \Pic (X)\\
@AAA    @AAA @AAA\\
H^2(\pi_1(Y,\eta),\hat \Bbb Z(1)) @>>> H^2(Y,\hat \Bbb Z(1)) @<<< \Pic (Y)\\
\endCD
$$
where the right and middle vertical maps are induced by the scheme morphism $X\to Y$ if $X$ is a formal $p$-adic germ, and the rigid morphism
$X\to Y^{\rig}$ and the comparison theorems between \'etale cohomology and rigid analytic 
\'etale cohomology in case $X$ is affinoid (cf. [Hansen], Theorem 1.8 and Theorem 1.9).  
The right horizontal maps are the cycle class maps arising from the Kummer exact sequence in \'etale topology
and the left lower horizontal map is an isomorphism arising from the Cartan-Leray spectral sequence (cf. [Mochizuki], Proposition 1.1). 
The pull-back of the class $[\pi_1(V,\eta)^{(\text c-\cn)}]\in  H^2(\pi_1(Y,\eta),\hat \Bbb Z(1))$ in $H^2(\pi_1(X,\eta),\hat \Bbb Z(1))$, 
where $V$ is the complement in $Y$ of the image of $S=\{x\}$ [cf. [Sa\"\i di1] 2.1.1 for the definition of $\pi_1(V,\eta)^{(\text c-\cn)}$], 
coincides with the class $[\pi_1(U,\eta)^{(\text c-\cn)}]$ (this follows from Lemma 2.1 and the various definitions). 
The class $[\pi_1(V,\eta)^{(\text c-\cn)}]\in  H^2(\pi_1(Y,\eta),\hat \Bbb Z(1))\isom H^2(Y,\hat \Bbb Z(1))$ coincides with the image of the Chern class of the line bundle 
$\Cal O(y)\in \Pic(Y)$ where $y\in Y(k)$ is the image of $x$ (cf. [Sa\"\i di1] proof of Lemma 2.3.1). Thus, the image of 
$[\pi_1(U,\eta)^{(\text c-\cn)}]$ in $H^2(X,\hat \Bbb Z(1))$ 
coincides with the image of the line bundle $\Cal O(x)\in \Pic(X)$ via the cycle class map $\Pic (X)\rightarrow H^2(X,\hat \Bbb Z(1))$ as claimed.

The Picard group $\Pic(X)$ is finite (cf. Theorem 4.1 in this paper in case $X$ is affinoid, and [Sa\"\i di2] Proposition 5.4 in case $X$ is a formal $p$-adic germ). 
In particular, the image of $[\pi_1(U,\eta)^{(\text c-\cn)}]$ in $H^2(X,\hat \Bbb Z(1))$
and hence the class $[\pi_1(U,\eta)^{(\text {c}-\cn)}]$ is a torsion element
of $H^2(\pi_1(X,\eta),\hat \Bbb Z(1))$. The class $[E_U]\in H^2(G_k,\hat \Bbb Z(1))$ of the group extension $E_U$ is the image of $[\pi_1(U,\eta)^{(\text {c}-\cn)}]$ under the retraction map 
$H^2(\pi_1(X,\eta),\hat \Bbb Z(1))@>s^{\star}>>H^2(G_k,\hat \Bbb Z(1))\isom \hat \Bbb Z$ induced by $s$. Hence the class $[E_U]$ is trivial since $\hat \Bbb Z$ is torsion free, and the group extension $E_U$ splits.
\qed
\enddemo

\proclaim{Theorem 2.3} {\bf (Lifting of sections to cuspidally abelian arithmetic fundamental groups)}. 
We use the above notations. There exists a section $s_U^{\ab}:G_k\to \pi_1(U,\eta)^{(\text c-\ab)}$ of the projection $\pi_1(U,\eta)^{(\text c-\ab)}\twoheadrightarrow G_k$ 
which {\bf lifts} the section $s$, i.e., which inserts in the following commutative diagram
$$
\CD
G_k @>s_U^{\text c-\ab} >> \pi_1(U,\eta)^{(\text c-\ab)}\\
@| @VVV\\
G_k @>s>> \pi_1(X,\eta)\\
\endCD
$$
where the right vertical map is the natural projection $ \pi_1(U,\eta)^{(\text c-\ab)}\twoheadrightarrow  \pi_1(X,\eta)$.
In particular, the set of sections of the 
projection $\pi_1(U,\eta)^{(\text c-\ab)}\twoheadrightarrow G_k$ which lift the section $s$ is non-empty, and is (up to conjugation by elements of $\Delta_U^{\ab}$) 
a torsor under $H^1(G_k,\Delta_U^{\ab})$.
\endproclaim

\demo{Proof} Let $\{H\}_{i\in I}$ be a projective system of open subgroups of $\pi_1(X,\eta)$ containing $s(G_k)$ such that $s(G_k)=\bigcap _{i\in I}H_i$. Thus, for $i\in I$, the open subgroup $H_i$ 
corresponds to an \'etale finite cover $X_i\to X$ with $X_i$ geometrically connected and $H_i$ is identified  with $\pi_1(X_i,\eta)$ which sits in the exact sequence 
$1\to \pi_1(X_{i},\eta)^{\geo}\to \pi_1(X_i,\eta)\to G_k\to 1$ (the geometric point; denote also $\eta$, of $X_i$ is induced by the geometric point $\eta$ of $X$). Further, the section $s$ induces a section $s_i:G_k\to \pi_1(X_i,\eta)$ of the projection $\pi_1(X_i,\eta)\twoheadrightarrow G_k$.
Let $U_i\defeq U\times _XX_i$ and $\pi_1(U_i,\eta)^{(\text c-\cn)}$ the (geometrically) cuspidally central arithmetic fundamental 
group of $U_i$ which sits in the exact sequence $1\to \Delta_{U_i}^{\cn}\to \pi_1(U_i,\eta)^{(\text c-\cn)}\to \pi_1(X_i,\eta)\to 1$.

Consider the following commutative diagrams
$$
\CD
1 @>>> \Delta_U^{\ab} @>>> \Cal E_U  @>>> G_k @>>> 1\\
@.         @VVV  @VVV  @VsVV \\
1 @>>> \Delta_U^{\ab} @>>> \pi_1(U,\eta)^{(\text c-\ab)} @>>> \pi_1(X,\eta) @>>> 1\\
\endCD
$$

and for $i\in I$
$$
\CD
1 @>>> \Delta_{U_i}^{\cn} @>>> E_{U_i}  @>>> G_k @>>> 1\\
@.         @VVV  @VVV  @Vs_iVV \\
1 @>>> \Delta_{U_i}^{\cn} @>>> \pi_1(U_i,\eta)^{(\text c-\cn)}@>>> \pi_1(X_i,\eta) @>>> 1\\
\endCD
$$
where the right squares are cartesian. Thus, $\Cal E_U$ (resp. $E_{U_i}$) is the pull-back of the group extension $\pi_1(U,\eta)^{(\text c-\ab)}$ [resp. $\pi_1(U_i,\eta)^{(\text c-\cn)}$] 
by the section $s$ (resp. $s_i$). 
There is a natural isomorphism $\Delta_U^{\ab}=\varprojlim _{i\in I}\Delta_{U_i}^{\cn}$ as follows from the facts that $\Delta_U=\Delta_{U_i}$, $\forall i\in I$, and given a finite quotient $\Delta_U^{\ab}\twoheadrightarrow H$ there exists $i\in I$ such that $\pi_1(X_i,\eta)^{\geo}$ acts trivially on $H$. 
Further, there is a natural isomorphism $\Cal E_U\isom \varprojlim _{i\in I} E_{U_i}$ (the transition maps in the projective limit being surjective). 
The existence of a section $s_U^{\text c-\ab}:G_k\to \pi_1(U,\eta)^{(\text c-\ab)}$ of the projection $\pi_1(U,\eta)^{(\text c-\ab)}\twoheadrightarrow G_k$ 
which lifts the section $s$ is equivalent to the splitting of the group extension $\Cal E_U$, and the set of those (possible) liftings $s_U^{\text c-\ab}$ is in one-to-one correspondence with the set of sections of the projection $\Cal E_U\twoheadrightarrow G_k$.
The natural projection $E_{U_i}\twoheadrightarrow G_k$ splits for all $i\in I$ (see proof of Proposition 2.2). We show the group extension $\Cal E_U$ splits.

Let $(P_j)_{j\in J}$ be a projective system of quotients $\Cal E_U
\twoheadrightarrow P_j$, where $P_j$ sits in an exact sequence $1\to F_j\to P_j\to G_k\to 1$ with $F_j$ finite, and $\Cal E_U=\varprojlim _{j\in J}P_j$. 
[More precisely, write $\Cal E_U$ as a projective limit of finite groups $\{\tilde P_j\}_{j\in J}$ where $\tilde P_j$ sits in an exact sequence 
$1\to F_j\to \tilde P_j\to G_j\to 1$ with $G_j$ a quotient of $G_k$ and $F_j$ a quotient of $\Ker(\Cal E_U\twoheadrightarrow G_k)$. 
Let $1\to F_j\to P_j\to G_k\to 1$ be the pull-back of the group extension $1\to F_j\to \tilde P_j\to G_j\to 1$ 
by $G_k\twoheadrightarrow G_j$. Then $\Cal E_U=\varprojlim _{j\in J}P_j$]. 
The set $\Sect(G_k,\Cal E_U)$ of group-theoretic sections of the projection $\Cal E_U\twoheadrightarrow G_k$ is naturally identified with 
the projective limit $\varprojlim _{j\in J}\Sect(G_k,P_j)$ 
of the sets $\Sect(G_k,P_j)$ of group-theoretic sections of the projections $P_j\twoheadrightarrow G_k$, $j\in J$. The set $\Sect(G_k,P_j)$ is non-empty, $\forall j\in J$. Indeed,
$P_j$ (being a quotient of $\Cal E_U$) is a quotient of $E_{U_i}$ for some $i\in I$, this quotient $E_{U_i}\twoheadrightarrow P_j$ commutes with the projections onto $G_k$, and we know the projection $E_{U_i}\twoheadrightarrow G_k$ splits, 
hence the projection $P_j\twoheadrightarrow G_k$ splits.
Moreover, the set $\Sect (G_k,P_j)$ is, up to conjugation by the elements of $F_j$, a torsor under the group $H^1(G_k,F_j)$ which is finite since $k$ is a $p$-adic local field [cf. [Neukirch-Schmidt-Wingberg], (7.1.8)Theorem(iii)]. Thus, 
 $\Sect (G_k,P_j)$ is a non-empty finite set. The set $\Sect (G_k,\Cal E_U)$ is non-empty being the projective limit of non-empty finite sets. This finishes the proof of Theorem 2.3.
\qed
\enddemo

Next let 
$$G_X^{(\text c-\ab)}\defeq G_X/\Ker (\Cal H_X\twoheadrightarrow \Cal H_X^{\ab})$$
(cf. exact sequence (1.4) for the definition of $\Cal H_X$). 
Thus, $G_X^{(\text c-\ab)}=\varprojlim _U \pi_1(U,\eta)^{(\text c-\ab)}$ where $U$ runs over all 
sub-schemes of $X$ as in  Theorem 2.3.

\proclaim{Theorem 2.4} {\bf (Lifting of sections to cuspidally abelian Galois groups)}. 
We use the above notations. There exists a section $s^{\text c-\ab}:G_k\to G_X^{(\text c-\ab)}$ of the projection $G_X^{(\text c-\ab)}\twoheadrightarrow G_k$ 
which {\bf lifts} the section $s$, i.e., which inserts in the following commutative diagram
$$
\CD
G_k @>s^{\text c-\ab} >> G_X^{(\text c-\ab)}\\
@| @VVV\\
G_k @>s>> \pi_1(X,\eta)\\
\endCD
$$
where the right vertical map is the natural projection $G_X^{(\text c-\ab)}\twoheadrightarrow \pi_1(X,\eta)$.
In particular, the set of sections of the 
projection $G_X^{(\text c-\ab)}\twoheadrightarrow G_k$ which lift the section $s$ is non-empty, and is (up to conjugation by elements of $\Cal H_X^{\ab}$) 
a torsor under $H^1(G_k,\Cal H_X^{\ab})$.
\endproclaim

\demo{Proof} The proof follows, using the natural identification $G_X^{\text c-\ab}\isom \varprojlim _{U} \pi_1(U,\eta)^{\text c-\ab}$ [where $U$ runs over all sub-schemes of $X$ as in  Theorem 2.3],
from Theorem 2.3 and a similar argument in our context to the one used in the proof of Theorem 2.3.5 in [Sa\"\i di1]. Alternatively, one can use Theorem 2.3 and a similar argument to the one used at the end of the proof of Theorem 2.3. 
\qed
\enddemo

The following is one of our main results in this section.

\proclaim {Theorem 2.5} 
Assume that $X$ admits a {\bf $k$-compatctification} $Y$ (cf. Notations).
If the projection $\pi_1(X,\eta)\twoheadrightarrow G_k$ {\bf splits} then $\index (Y)=1$. 
\endproclaim

\demo{Proof}  Assume that the projection $\pi_1(X,\eta)\twoheadrightarrow G_k$ splits and let $s:G_k\to \pi_1(X,\eta)$ be a section of this projection.
By Theorem 2.4 there exists a section $s^{\text c-\ab}:G_k\to G_X^{(\text c-\ab)}$ of the projection $G_X^{(\text c-\ab)}\twoheadrightarrow G_k$ 
which lifts the section $s$. The section $s^{\text c-\ab}$ induces naturally a section $\tilde s:G_k\to G_X^{(\ab)}$ of the projection $G_X^{(\ab)}\twoheadrightarrow G_k$ 
(see $\S1$ for the definition of $G_X^{(\ab)}$ and note that $G_X^{(\ab)}$ is a quotient of $G_X^{(\text c-\ab)}$). 
Let $G_Y\defeq \Gal (\overline K/K)$ be the absolute Galois group of the function field $K$ of $Y$ and 
$G_Y^{(\ab)}\defeq G_Y/\Ker [ \Gal (\overline K/K\bar k)\twoheadrightarrow \Gal (\overline K/K\bar k)^{\ab} ]$ 
its geometrically abelian quotient. We have a commutative diagram
$$
\CD
G_X^{(\ab)} @>>> G_k\\
@VVV @|\\
G_Y^{(\ab)}@>>> G_k\\
\endCD
$$
where the left vertical map is induced by the natural map $G_X\to G_Y$, which is induced by the scheme morphism $X\to Y$ in case $X$ is a formal $p$-adic germ, and
by the rigid analytic morphism $X\to Y^{\rig}$ and the rigid GAGA functor in case $X$ is affinoid.
The section $\tilde s:G_k\to G_X^{(\ab)}$ induces a section $s^{\dag}:G_k\to G_Y^{(\ab)}$ of the projection $G_Y^{(\ab)}\twoheadrightarrow G_k$ (cf. above diagram). 
The existence of the section $s^{\dag}$ implies that $\index (Y)=1$ as was observed by Esnault and Wittenberg (see [Esnault-Wittenberg] Remark 2.3(ii), and 
[Sa\"\i di5] Theorem A for a more general result). 
\qed
\enddemo

\subhead
\S3. Geometric sections of arithmetic fundamental groups 
\endsubhead
We investigate {\bf geometric} sections of the projection $\pi_1(X,\eta)\twoheadrightarrow G_k$ (relative to a fixed compactification of $X$).
We use the notations in $\S0$, $\S1$, $\S2$. We further assume that $X$ possesses a {\bf $k$-compactification} $Y$ with $Y$ {\bf hyperbolic} (cf. Notations) which is fixed throughout $\S3$.

Let 
$$s:G_k\to  \pi_1(X,\eta)$$ 
be a {\bf section} of the projection $\pi_1(X,\eta)\twoheadrightarrow  G_k$ fixed throughout $\S3$, which induces a {\bf (local) section}
$$s_Y:G_k\to  \pi_1(Y,\eta)$$ 
of the projection $\pi_1(Y,\eta)\twoheadrightarrow  G_k$ [cf. diagram (0.1) and $\S0$]. 

We have an exact sequence 
$$1\to \Cal I_Y\to G_Y\to \pi_1(Y,\eta)\to 1$$ 
where $G_Y=\Gal (\overline K/K)$ is the absolute Galois group of the function field $K$ of $Y$ and $\Cal I_Y\defeq \Ker[G_Y\twoheadrightarrow \pi_1(Y,\eta)]$.
Let 
$$G_Y^{(\text c-\ab)}\defeq G_Y/\Ker (\Cal I_Y\twoheadrightarrow \Cal I_Y^{\ab}).$$ 
Thus $G_Y^{(\text c-\ab)}=\varprojlim _V \pi_1(V,\eta)^{(\text c-\ab)}$ where $V$ runs over all 
open sub-schemes of $Y$ [cf. [Sa\"\i di1], 2.1.1, for the definition of $\pi_1(V,\eta)^{(\text c-\ab)}$].

\proclaim{Theorem 3.1} {\bf (Lifting of sections to cuspidally abelian Galois groups)}.
We use the above notations. The followings hold.

(i)\ There exists a section $s_Y^{\text c-\ab}:G_k\to G_Y^{(\text c-\ab)}$ of the projection $G_Y^{(\text c-\ab)}\twoheadrightarrow G_k$ 
which {\bf lifts} the section $s_Y:G_k\to \pi_1(Y,\eta)$, i.e., which inserts in the following commutative diagram
$$
\CD
G_k @>s_Y^{\text c-\ab} >> G_Y^{(\text c-\ab)}\\
@| @VVV\\
G_k @>s_Y>> \pi_1(Y,\eta)\\
\endCD
$$
where the right vertical map is the natural projection
$G_Y^{(\text c-\ab)}\twoheadrightarrow \pi_1(Y,\eta)$. In particular, the set of sections of the 
projection $G_Y^{(\text c-\ab)}\twoheadrightarrow G_k$ which lift the section $s_Y$ is non-empty, and is (up to conjugation by elements of $\Cal I_Y^{\ab}$) 
a torsor under $H^1(G_k,\Cal I_Y^{\ab})$. 

(ii)\ The (local) section $s_Y:G_k\to \pi_1(Y,\eta)$ is {\bf uniformly orthogonal to $\Pic$} in the sense of [Sa\"\i di1], Definition 1.4.1.
\endproclaim

\demo{Proof} Assertion (i) follows from Theorem 2.4 and the fact that there exists a natural homomorphism $G_X^{(\text c-\ab)}\to G_Y^{(\text c-\ab)}$, induced by the natural homomorphism
$G_X\to G_Y$, which commutes with the projections to $G_k$. Assertion (ii) follows from assertion (i) and Theorem 2.3.5 in [Sa\"\i di1]. 
\qed
\enddemo



Consider the following push-out diagram
$$
\CD
1  @>>> \Cal H_X    @>>> G_X    @>>>   \pi_1(X,\eta) @>>> 1\\
@.        @VVV           @VVV        @| \\
1  @>>> \Cal H_{X,1/p^2}    @>>> G_X ^{(1/p^2-\sol)}    @>>>   \pi_1(X,\eta) @>>> 1\\
\endCD
$$
where $\Cal H_{X,1/p^2}$ is the {\bf maximal $1/p^2$-th solvable quotient} of $\Cal H_X$ and $G_X ^{(1/p^2-\sol)}\defeq G_X/\Ker ( \Cal H_X \twoheadrightarrow \Cal H_{X,1/p^2})$.
Thus, $\Cal H_{X,1/p^2}$ is the maximal quotient of $\Cal H_X$ which is abelian and annihilated by $p^2$
 (cf. [Sa\"\i di2], 1.2, for more details). 
We have a commutative diagram of exact sequences

$$
\CD
1  @>>> \Cal H_{X,1/p^2}    @>>> G_X ^{(1/p^2-\sol)}    @>>>   \pi_1(X,\eta) @>>> 1\\
@.        @VVV           @VVV        @VVV \\
1  @>>> \Cal I_{Y,1/p^2}    @>>> G_Y ^{(1/p^2-\sol)}    @>>>   \pi_1(Y,\eta) @>>> 1\\
\endCD
\tag 3.1
$$
which is induced by the natural homomorphism $G_X\to G_Y$, where $G_Y ^{(1/p^2-\sol)}$ is defined in a similar way to $G_X ^{(1/p^2-\sol)}$.  
More precisely,  $\Cal I_{Y,1/p^2}$ is the maximal quotient of $\Cal I_Y$ which is abelian and annihilated by $p^2$ and 
$G_Y ^{(1/p^2-\sol)}\defeq G_Y/\Ker ( \Cal I_Y \twoheadrightarrow \Cal I_{Y,1/p^2})$ is the {\bf geometrically cuspidally $1/p^2$-th step solvable quotient}
of $G_Y$ [cf. [Sa\"\i di2], 3.1, recall the exact sequence $1 \to  \Cal I_{Y}   \to G_Y   \to   \pi_1(Y,\eta) \to  1$].

The following Proposition 3.2, item (i), is weaker than (and follows from) Theorem 2.4, we state it in connection with Theorem  3.5.2  in this section.


\proclaim {Proposition 3.2 (Lifting of sections to cuspidally $1/p^2$-th step solvable Galois groups)} 
We use the above notations. The followings hold.

(i) There exists a section $\tilde s:G_k\to G_X^{(1/p^2-\sol)}$ of the projection $G_X^{(1/p^2-\sol)}\twoheadrightarrow G_k$ 
which {\bf lifts} the section $s:G_k\to \pi_1(X,\eta)$, i.e., which inserts in the following commutative diagram
$$
\CD
G_k @>\tilde s>> G_X^{(1/p^2-\sol)}\\
@| @VVV\\
G_k @>s>> \pi_1(X,\eta)\\
\endCD
$$
where the right vertical map is the natural projection
$G_X^{(1/p^2-\sol)}\twoheadrightarrow \pi_1(X,\eta)$.
In particular, the set of sections of the projection $G_X^{(1/p^2-\sol)}\twoheadrightarrow G_k$ which lift the section $s$ is non-empty, and is (up to conjugation by elements of $\Cal H_{X,1/p^2}$) 
a torsor under $H^1(G_k,\Cal H_{X,1/p^2})$.

(ii) The section $\tilde s:G_k\to G_X^{(1/p^2-\sol)}$ in (i) induces a section $\tilde s_Y:G_k\to G_Y^{(1/p^2-\sol)}$ of the projection 
$G_Y^{(1/p^2-\sol)}\twoheadrightarrow G_k$ which {\bf lifts} the section $s_Y:G_k\to \pi_1(Y,\eta)$. In particular, the (local) section $s_Y:G_k\to \pi_1(Y,\eta)$ 
is {\bf uniformly orthogonal to $\Pic$ mod-$p^2$} in the sense of [Sa\"\i di2], Definition 3.4.1.
\endproclaim

\demo{Proof} Assertion (i) follows from Theorem 2.4 and the fact that there exists a natural projection $G_X^{(\text c-\ab)}\twoheadrightarrow  G_X^{(1/p^2-\sol)}$ which commutes with the projections onto $G_k$. Assertion (ii) follows from (i) and the fact that there exists a natural homomorphism $G_X^{(1/p^2-\sol)}\to G_Y^{(1/p^2-\sol)}$, induced by the 
homomorphism $G_X\to G_Y$, which commutes with the projections onto $G_k$ (cf. diagram (3.1), and [Sa\"\i di2] Theorem 3.4.4).
\qed
\enddemo

\subhead 
3.3
\endsubhead
Write   
$$\Pi_Y[X]\defeq \varprojlim _{T\subset Y\setminus X}\pi_1(Y\setminus T,\eta)$$
and
$$\Pi_Y[X]^{\geo}\defeq \varprojlim _{T\subset Y\setminus X}\pi_1(Y\setminus T,\eta)^{\geo}$$
where the limits are over all subsets $T$ consisting of finitely many closed points of $Y$ {\bf not in} $X$ (cf. Notations),
$Y\setminus T$ is the corresponding (affine if $T$ is non-empty) curve,
and  $\pi_1(Y\setminus T,\eta)^{\geo}\defeq \Ker [\pi_1(Y\setminus T,\eta)\twoheadrightarrow G_k]$.
We have the following commutative diagram of exact sequences
$$
\CD
1@>>>\pi_1(X,\eta)^{\geo}@>>> \pi_1(X,\eta)@>>> G_k @>>> 1 \\
@.@VVV  @VVV @|\\
1 @>>>     \Pi_Y[X]^{\geo} @>>>  \Pi_Y[X] @>>> G_k@>>> 1\\
@.@VVV  @VVV @|\\
1@>>> \pi_1(Y_{\bar k},\bar \eta)@>>> \pi_1(Y,\eta)@>>> G_k@>>> 1\\
\endCD
\tag 3.2
$$
where the middle upper map is induced by the rigid analytic morphism $X\to Y^{\rig}$ and the rigid GAGA functor 
in case $X$ is {\bf affinoid}, and the scheme morphism $X\to Y$ in case $X$ is a {\bf formal $p$-adic germ}.
The left and middle lower vertical maps are the natural projections (they are surjective).

\proclaim{Proposition 3.3.1} We use the above notations. The left and middle upper vertical maps in diagram (3.2) are {\bf injective} in the case $X$ is {\bf affinoid}.
\endproclaim

\demo{Proof} The first assertion follows from Theorem A in [Sa\"\i di3] (see the comments in the proof of Proposition 1.1). 
The second assertion follows from the first and the commutativity of the upper part in diagram (3.2).
\qed
\enddemo

The section $s:G_k\to  \pi_1(X,\eta)$ induces a section (denoted also $s$)
$$s:G_k\to   \Pi_Y[X]$$ 
of the projections $\Pi_Y[X] \twoheadrightarrow G_k$ 
[cf. diagram (3.2)].

\definition {Definition 3.3.2} We say that the section $s$ is {\bf geometric}, relative to $Y$, if the image $s(G_k)$ of the section
$s:G_k\to   \Pi_Y[X]$ is contained in a decomposition group 
$D_x\subset \Pi_Y[X]$ associated to a {\bf rational} point $x\in Y(k)$.
\enddefinition

Note that if $s$ is geometric in the above sense, associated to $x\in Y(k)$, then the (local) section $s_Y:G_k\to   \pi_1(Y,\eta)$ 
of the projection $\pi_1(Y,\eta)\twoheadrightarrow G_k$ induced by $s$ is geometric and is associated to $x\in Y(k)$, i.e.,
$s_Y(G_k)$ is contained in (hence equal to) a decomposition group $D_x\subset \pi_1(Y,\eta)$ associated to $x$.

\subhead 3.4
\endsubhead
In this sub-section we assume that $X=\Spec (A\otimes _{\Cal O_k}k)$ is a {\bf formal $p$-adic germ}. 

Let $\Cal Y\to \Spec \Cal O_k$ be a model of $Y$, $y\in \Cal Y^{\cl}$ a closed point, and $\hat \Cal O_{\Cal Y,y}\isom A$ an isomorphism (cf. loc. cit.). 
Let $\Cal Y_F\defeq \Cal Y\times _{\Spec \Cal O_k}\Spec F$ be the special fibre of $\Cal Y$. Consider the following assumption {\bf (*)}:

\bigskip
{\bf (*)} {\it The gcd of the total multiplicities of the irreducible components of $\Cal Y_F$ is $1$}.

\bigskip
Let $\xi$ be a geometric point of $\Cal Y_F$ with values in the generic point of an 
irreducible component $Y_{i_0}$ of $\Cal Y_F$. Thus $\xi$ determines an algebraic closure $\overline F$ of $F$.
We have the following commutative diagram of exact sequences
$$
\CD
1@>>>\pi_1(X,\eta)^{\geo}@>>> \pi_1(X,\eta)@>>> G_k @>>> 1 \\
@.@VVV  @VVV @|\\
1@>>> \pi_1(Y_{\bar k},\bar \eta)@>>> \pi_1(Y,\eta)@>>> G_k@>>> 1\\
@.@VVV  @VVV @VVV\\
1@>>> \pi_1(\Cal Y_{\overline F},\bar \xi)@>>> \pi_1(\Cal Y_F,\xi)@>>> G_F@>>> 1\\
\endCD
\tag 3.3
$$
where the middle upper map is induced by the scheme morphism $X\to Y$, the lower middle map 
(which is defined up to conjugation) is a specialisation map, 
$\pi_1(\Cal Y_F,\xi)$ [resp. $\pi_1(\Cal Y_{\overline F},\bar \xi)$] is the fundamental group of $\Cal Y$
(resp. $\Cal Y_{\overline F}\defeq \Cal Y\times _{\Spec \Cal O_k}\Spec \overline F$) with geometric point $\xi$ (resp. $\bar \xi$ which is induced by $\xi$),
$G_F\defeq \Gal (\overline F/F)$, and the lower right vertical map is the natural projection $G_k\twoheadrightarrow G_F$ (cf. [Sa\"\i di6], diagram (0.1), and the discussion thereafter).
The left (hence also the middle) lower vertical map in diagram (3.3) is surjective under the assumption {\bf (*)} (cf. loc. cit. and the references therein).

The section $s:G_k\to \pi_1(X,\eta)$ induces  the (local) section $s_Y:G_k\to \pi_1(Y,\eta)$ 
of the projection $\pi_1(Y,\eta)\twoheadrightarrow G_k$, as well as a homomorphism 
$$\tilde s:G_k\to \pi_1(\Cal Y_F,\xi)$$ 
obtained by composing the section 
$s_Y:G_k\to \pi_1(Y,\eta)$ with the specialisation map
$\pi_1(Y,\eta)\twoheadrightarrow \pi_1(\Cal Y_F,\xi)$ in diagram (3.3).

\proclaim {Lemma 3.4.1} We use the above notations. The followings hold.

\noindent
(i) The closed point $y\in \Cal Y^{\cl}$ is an $F$-{\bf rational} point.

\noindent
(ii) The section $s_Y$ is {\bf unramified}: the homomorphism $\tilde s:G_k\to \pi_1(\Cal Y_F,\xi)$ 
factors through $G_F$ and induces a section $\bar s_Y:G_F\to \pi_1(\Cal Y_F,\xi)$ of the natural projection
$\pi_1(\Cal Y_F,\xi)\twoheadrightarrow  G_F$.

\noindent
(iii) The section $\bar s_Y:G_F\to \pi_1(\Cal Y_F,\xi)$ in (ii) is {\bf geometric} and arises from the rational point $y$, i.e., arises from the scheme-theoretic morphism
$y:\Spec F\to \Cal Y_F$.

\noindent
(iv) Assume that $\Cal Y$ is regular. Then condition {\bf (*)} holds.

\endproclaim

\demo{Proof} Assertion (i) is clear (recall $\hat \Cal O_{\Cal Y,y}\isom A$), it also follows from (ii).
We prove (ii). 

We have a commutative diagram of scheme morphisms
$$
\CD
X  @>>> Y\\
@VVV   @VVV\\
\Spec A @>>> \Cal Y\\
@AAA    @AAA\\
\Spec (F)@>y>> \Cal Y_F\\
\endCD
\tag 3.4
$$
where the lower horizontal morphism is induced by the closed point $y$ of $\Cal Y_F$,
and the lower vertical morphisms are closed immersions. This diagram gives rise to a commutative diagram of 
homomorphisms between fundamental groups

$$
\CD
\pi_1(X,\eta)  @>>> \pi_1(Y,\eta)\\
@VVV     @VVV\\
\pi_1(\Spec A,\eta) @>>> \pi_1(\Cal Y,\eta)\\
@A\tau _yAA   @A{\sigma}AA\\
G_{F} @>s_y>> \pi_1(\Cal Y_F,\xi)\\
\endCD
\tag 3.5
$$
where the lower horizontal map is a section
of the projection $\pi_1(\Cal Y_F,\xi)\twoheadrightarrow G_F$ arising from the $F$-rational point $y\in \Cal Y_F$, and is defined up to conjugation,
the lower vertical maps are induced by the lower vertical maps in diagram (3.4) (they are defined up to conjugation) and are isomorphisms (cf. [Grothendieck], Expos\'e X, Th\'eor\`eme 2.1,
for the right vertical map $\sigma$ being an isomorphism).
 Further, the composite $\psi:\pi_1(X,\eta)\to \pi_1(\Spec A,\eta) @>\tau_{y}^{-1}>> G_{F} @>s_y>>
\pi_1(\Cal Y_F,\xi)$ is the composite of the middle vertical maps in diagram (3.3) as follows from the definition of the specialisation map $\pi_1(Y,\eta)\to \pi_1(\Cal Y_F,\xi)$:
this map is the composite of the maps $\pi_1(Y,\eta)\to \pi_1(\Cal Y,\eta)@>{\sigma^{-1}}>>\pi_1(\Cal Y_F,\xi)$.
In particular, the homomorphism $\tilde s:G_k\to \pi_1(\Cal Y_F,\xi)$ 
factors through $G_F$ and induces a section $\bar s_Y:G_F\to \pi_1(\Cal Y_F,\xi)$ of the natural projection $\pi_1(\Cal Y_F,\xi)\twoheadrightarrow  G_F$. This shows (ii).
The section $\bar s_Y$ coincides (up to conjugation) with the section $G_{F} @>s_y>> \pi_1(\Cal Y_F,\xi)$ in diagram (3.5), hence is geometric and arises from the $F$-rational point $y$ as claimed in (iii).
The last assertion follows from Theorem 2.5 and the well known fact that if $\Cal Y$ is regular then the gcd of the total multiplicities of the irreducible components of $\Cal Y_F$
divides $\index (Y)$ (cf. for example [Gabber-Liu-Lorenzini], Theorem 8.2 and Remark 8.6). 
\qed
\enddemo

\definition {Remark 3.4.2}
Assume that the morphism $\Cal Y\to \Spec \Cal O_k$ is {\bf smooth}. If $s$ is geometric, and arises from the rational point $x\in Y(k)$ (cf. Definition 3.3.2), it follows from Lemma 3.4.1(iii) 
and the fact that $\Cal Y_F$ is hyperbolic that the point $x$ specialises in $y$ necessarily [cf. [Tamagawa], Proposition (2.8)(i)]. In particular, the point $x$ is the image of a (unique) $k$-rational point $\tilde x\in X(k)$ via the morphism $X\to Y$. The fact that $s_Y(G_k)=D_x\subset \pi_1(Y,\eta)$ doesn't imply a priori that the image $s(G_k)$ via the section 
$s:G_k\to \pi_1(X,\eta)$ is contained in a decomposition group  $D_{\tilde x}\subset \pi_1(X,\eta)$ associated to $\tilde x$.
\enddefinition

\subhead {3.5}
\endsubhead
Let $H\subset \Pi_Y[X]$ be an open subgroup with $s(G_k)\subset H$ [recall $s:G_k\to \Pi_Y[X]$ is the section induced by $s:G_k\to \pi_1(X,\eta)$].  
Thus, $H$ corresponds to a (possibly ramified) finite cover $Y'\to Y$ with $Y'$ geometrically connected. 
Let $H' \subset \pi_1(X,\eta)$ be the inverse image of $H$ via the homomorphism $\pi_1(X,\eta)\to \Pi_Y[X]$ (cf. diagram (3.2)). 
Thus, $H'$ is an open subgroup of  $\pi_1(X,\eta)$ containing the image of the section $s:G_k\to \pi_1(X,\eta)$ and corresponds to an \'etale cover 
$X'\to X$ with $X'$ geometrically connected. There is a natural morphism $X'\to (Y')^{\rig}$ of rigid analytic spaces in case $X$ is {\bf affinoid}, and a natural scheme morphism 
$X'\to Y'$ in case $X$ is a {\bf formal $p$-adic germ}. The generic point $\eta$ induces naturally a generic point (denoted also $\eta$) of $X'$ and $Y'$. Further, we have a natural identification $H'=\pi_1(X',\eta)$ and a natural homomorphism $\pi_1(X',\eta)\to \pi_1(Y',\eta)$ which commutes with the projections onto $G_k$. 

The section $s:G_k\to \pi_1(X,\eta)$ induces naturally sections $s':G_k\to \pi_1(X',\eta)$ and $s_{Y'}:G_k\to \pi_1(Y',\eta)$ of the natural projections  $\pi_1(X',\eta)\twoheadrightarrow G_k$
and $\pi_1(Y',\eta)\twoheadrightarrow G_k$; respectively. The section $s':G_k\to \pi_1(X',\eta)$ lifts to a section $\tilde s':G_k\to G_{X'}^{(1/p^2-\sol)}$ of the projection $G_{X'}^{(1/p^2-\sol)}\twoheadrightarrow G_k$ and induces a section $\tilde s_{Y'}:G_k\to G_{Y'}^{(1/p^2-\sol)}$ of the projection $G_{Y'}^{(1/p^2-\sol)}\twoheadrightarrow G_k$ (cf. Proposition 3.2). 
Let $F\subset G_{Y'}^{(1/p^2-\sol)}$ be an open subgroup with $\tilde s_{Y'}(G_k)\subset F$. Thus $F$ corresponds to a (possibly ramified) finite cover $Y''\to Y'$ with $Y''$ geometrically connected. The generic point $\eta$ induces naturally a generic point (denoted also $\eta$) of $Y''$. Write $\pi_1(Y'',\eta)^{(1/p-\sol)}$ for 
{\bf the geometrically $1/p$-th step solvable quotient} of $\pi_1(Y'',\eta)$ which sits in the following exact sequence
$$1@>>> \pi_1(Y''_{\bar k},\bar \eta)_{1/p}@>>> \pi_1(Y'',\eta)^{(1/p-\sol)} @>>> G_k@>>> 1,\tag 3.6$$
where $\pi_1(Y''_{\bar k},\bar \eta)_{1/p}$ is the maximal $1/p$-th step solvable quotient of $\pi_1(Y''_{\bar k},\bar \eta)$
(cf. [Sa\"\i di2], 1.2) and the generic point $\bar \eta$ is induced by $\eta$. Thus $\pi_1(Y''_{\bar k},\bar \eta)_{1/p}$ is the maximal quotient of $\pi_1(Y''_{\bar k},\bar \eta)$ which is abelian and annihilated by $p$ (cf. loc. cit.)

\definition {Definition 3.5.1} 
We use the above notations. We say that the section $s$ is {\bf admissible}, relative to $Y$, if for every open subgroup
$H\subset \Pi_Y[X]$ with $s(G_k)\subset H$, corresponding to (a possibly ramified) cover $Y'\to Y$, 
the following holds. There exists a section 
$\tilde s_{Y'}:G_k\to G_{Y'}^{(1/p^2-\sol)}$ of the projection $G_{Y'}^{(1/p^2-\sol)}\twoheadrightarrow G_k$ [such a section exists unconditionally (see above discussion)]
satisfying  the following property:
\smallskip
{\it For each open subgroup $F\subset G_{Y'}^{(1/p^2-\sol)}$ with $\tilde s_{Y'}(G_k)\subset F$, corresponding to a (possibly ramified) cover $Y''\to Y'$ with $Y''$ geometrically connected, the natural projection $\pi_1(Y'',\eta)^{(1/p-\sol)}\twoheadrightarrow G_k$ {\bf splits} (cf. above discussion)}.
\smallskip
Note that this latter condition is equivalent to (cf. [Sa\"\i di2] Lemma 3.4.8): 
\smallskip
{\it The class of $\Pic^1_{Y''}$ in $H^1(G_k,\Pic^0_{Y''})$ is divisible by $p$}.
\enddefinition

Our main result in this section is the following.

\proclaim {Theorem 3.5.2} We use the above notations. The section $s:G_k\to \pi_1(X,\eta)$ is {\bf geometric} relative to $Y$ (cf. Definition 3.3.2) {\bf if and only if} $s$  is {\bf admissible} relative to $Y$ (cf. Definition 3.5.1).
\endproclaim

\demo{Proof} Assume first that the section $s:G_k\to \pi_1(X,\eta)$ is admissible (relative to $Y$).
We prove that $s$ is geometric (relative to $Y$). Using a well-known limit argument due to Tamagawa [cf. [Tamagawa], Proposition 2.8(iv)] it suffices to show the following. 
For every open subgroup $H\subset \Pi_Y[X]$ with $s(G_k)\subset H$, corresponding to (a possibly ramified) cover $Y'\to Y$ with $Y'$ hyperbolic,
$Y'(k)\neq \emptyset $ holds. By assumption there exists a section $\tilde s_{Y'}:G_k\to G_{Y'}^{(1/p^2-\sol)}$ of the projection $G_{Y'}^{(1/p^2-\sol)}\twoheadrightarrow G_k$ satisfying the condition in Definition 3.5.1. In [Sa\"\i di2], 3.3, we defined a certain quotient $G_{Y'}\twoheadrightarrow G_{Y'}^{(p,2)}\twoheadrightarrow G_{Y'}^{(1/p^2-\sol)}$ of $G_{Y'}$
(we refer to loc. cit. for more details on the definition of $G_{Y'}^{(p,2)}$). Let $F\subset G_{Y'}^{(1/p^2-\sol)}$ be an open subgroup with $\tilde s_{Y'}(G_k)\subset F$ corresponding to a (possibly ramified) cover $Y''\to Y'$ with $Y''$ geometrically connected. 
By assumption the natural projection $\pi_1(Y'',\eta)^{(1/p-\sol)}\twoheadrightarrow G_k$ splits (cf. Definition 3.5.1). This latter condition (for every $F$ as above) implies that (in fact is equivalent to) the section $\tilde s_{Y'}:G_k\to G_{Y'}^{(1/p^2-\sol)}$ lifts to a section $s_{Y'}^{\dag}:G_k\to G_{Y'}^{(p,2)}$ of the projection $G_{Y'}^{(p,2)}\twoheadrightarrow G_k$ (cf. [Sa\"\i di2] Theorem 3.4.10 and Lemma 3.4.8). Further, the existence of the section $s_{Y'}^{\dag}:G_k\to G_{Y'}^{(p,2)}$ as above implies that $Y'(k)\neq \emptyset$ by [Sa\"\i di2], Proposition 4.6, as required.

Next, we assume that $s$ is geometric (relative to $Y$) and prove that $s$ is admissible (relative to $Y$). By assumption $s(G_k)$ is contained in $D_x\subset \Pi_Y[X]$ where $D_x$ is a decomposition group associated to a rational point $x\in Y(k)$. Let
$H\subset \Pi_Y[X]$ be an open subgroup with $s(G_k)\subset H$ corresponding to (a possibly ramified) cover $Y'\to Y$. Then $Y'(k)\neq \emptyset$. 
A rational point $x'\in Y'(k)$ gives rise to a section $\tilde s_{Y'}:G_k\to G_{Y'}^{(1/p^2-\sol)}$ of the projection $G_{Y'}^{(1/p^2-\sol)}\twoheadrightarrow G_k$.
Let $F\subset G_{Y'}^{(1/p^2-\sol)}$ be an open subgroup with $\tilde s_{Y'}(G_k)\subset F$ corresponding to a (possibly ramified) cover $Y''\to Y'$ with $Y''$ geometrically connected. 
Then $Y''(k)\neq \emptyset$ holds since the section  $\tilde s_{Y'}:G_k\to G_{Y'}^{(1/p^2-\sol)}$ 
arises from the rational point $x'$ and $\tilde s_{Y'}(G_k)\subset F$. In particular, the natural projection $\pi_1(Y'',\eta)\twoheadrightarrow G_k$, and a fortiori the projection $\pi_1(Y'',\eta)^{(1/p-\sol)}\twoheadrightarrow G_k$, splits. Thus $s$ is admissible as required. 
\qed
\enddemo

\subhead
\S 4. Picard groups of affinoid $p$-adic curves
\endsubhead
The following is our main result in this section, it may be of interest independently of the topics discussed in $\S1$, $\S2$, and $\S3$.

\proclaim{Proposition 4.1}  Let $X=\Sp(A)$ be a {\bf smooth} and geometrically connected  {\bf $k$-affinoid curve}. Then the Picard group $\Pic(X)$ is {\bf finite}.
\endproclaim

The rest of this section is devoted to the proof of Proposition 4.1.

Let $\Cal X=\Spf B$ be an excellent normal $\Cal O_k$-formal scheme of finite type
with generic fibre $X$, i.e., $A=B\otimes _R k$.
 Write $\Cal X^{\reg}$ for the set of regular points of $\Cal X$. Thus 
$\Cal X\setminus \Cal X^{\reg}=\{z_1,\ldots,z_t\}$ consists of finitely many closed points of $\Cal X$. By Lipman's theorem of resolution of singularities for excellent $2$-dimensional schemes there exists a birational and proper morphism $\lambda:\Cal S\to \Cal X$ with $\Cal S$ regular and $\lambda^{-1}(\Cal X^{\reg})\to \Cal X^{\reg}$ an isomorphism (cf. [Lipman]. Here we view $\Cal X$ as the ordinary affine scheme $\Spec B$). For $n\ge 1$, 
write $B_n\defeq B/(\pi^n)$, $\Cal X_n\defeq \Spec B_n$, and $\Cal S_n\defeq \Cal S\times _{\Cal X}\Cal X_n$.
Further denote $\Cal X_0\defeq \Cal X_n^{\red}$, 
and $\Cal S_0\defeq \Cal S_n^{\red}$. Thus $\Cal X_0$ and $\Cal S_0$ are one dimensional reduced schemes over $F$.
Further, there exists a morphism $\lambda:\Cal S\to \Cal X$ as above with $S_0$ a divisor with strict normal crossings (cf. [Cossart-Jannsen-Saito], Corollary 0.4), which we assume from now on.

We have a surjective homomorphism $\Pic(\Cal X^{\reg})\twoheadrightarrow \Pic(X)$. To prove $\Pic(X)$ is finite it suffices to prove that $\Pic(\Cal X^{\reg})$ is finite.
For each singular point $z_i$ of $\Cal X$ let $E_i\defeq \lambda^{-1}(z_i)^{\red}$ and $\{D_{i,j}\}_{1\le j\le n_i}$ the set of irreducible components  of $E_i$, $1\le i\le t$. 
Thus $E_i$ is a reduced proper curve over the residue field $k(z_i)$ at $z_i$ which is a finite field. 
We have an exact sequence
$$M\defeq \oplus_{i=1}^t(\oplus _{j=1}^{n_i}\Bbb Z)@>\beta>> \Pic(\Cal S)\to \Pic(\Cal X^{\reg})\to 0$$
where $\beta$ maps the copy of $\Bbb Z$ indexed by the pair $(i,j)$ to the class of the divisor $D_{i,j}$. Further we have an isomorphism 
$$\Pic (\Cal S)\isom \varprojlim _{n\ge 1}\Pic (\Cal S_n)$$
(cf. [EGA III], premi\`ere partie, Corollaire 5.1.6).

\proclaim {Lemma 4.2} We use notations as above. To prove that $\Pic(\Cal X^{\reg})$ is {\bf finite} it suffices to prove the following two assertions:

{\bf (A)}\ The cokernel of the composite map
$$\phi_n:M \defeq \oplus_{i=1}^t(\oplus _{j=1}^{n_i}\Bbb Z) @>\beta>> \Pic(\Cal S)\to \Pic (\Cal S_n)$$
is {\bf finite} for $n\ge 1$.

{\bf (B)}\ There exists $n_0>0$ such that the map
$$\Pic (\Cal S_{n+1})\to \Pic (\Cal S_n)$$
is an {\bf isomorphism} for $n>n_0$.
\endproclaim

\demo{Proof of Lemma 4.2} Follows from the above discussion and the fact that we have an exact sequence
$$M \to \varprojlim _{n\ge 1}\Pic (\Cal S_n) \to \varprojlim _{n\ge 1} \text {coker}(\phi_n)\to 0$$
where the first map is induced by the maps $\phi_n:M\to \Pic (\Cal S_n)$, $n\ge 1$, and $\varprojlim _{n\ge 1} \text {coker}(\phi_n)$ is finite if assertions {\bf (A)} and {\bf (B)} are satisfied. 

This finishes the proof of Lemma 4.2.
\qed
\enddemo

The rest of this section is devoted to the proofs of assertions {\bf (A)} and {\bf (B)}.

\demo{Proof of assertion {\bf (A)}} Let $\{\eta_r\}_{r=1}^s$ be the generic points of $\Cal X_0$,
$\rho:\Cal S_0^{\text {nor}}\to \Cal S_0$ the morphism of normalisation, $\widetilde E_i\defeq \rho^{-1}(E_i)$, $1\le i\le t$, and 
$H_r=\overline{\{\eta_r\}}$ the closure in $\Cal S_0^{\text {nor}}$ of the (inverse image in $\Cal S_0$ of the) generic point $\eta_r$ of $\Cal X_0$, $1\le r\le s$. 
Thus $H_r$ is a connected affine normal one dimensional scheme over $F$. 
Let 
$$d:\Pic (\Cal S_0)@>\rho^{*}>> \Pic (\Cal S_0^{\text {nor}})  @>\text {deg}>> M=\oplus_{i=1}^t(\oplus _{j=1}^{n_i}\Bbb Z)$$
be the composite map where the first map is the pullback of line bundles via the normalisation morphism 
$\rho:\Cal S_0^{\text {nor}}\to \Cal S_0$, and the 
map $\text {deg}$ is obtained by taking the degree of a line bundle on each irreducible component $D_{i,j}$ of $E_i$. 

\proclaim {Claim 1} $\text {ker}(d)$ is {\bf finite}
\endproclaim
\demo {Proof of Claim 1} We have a commutative diagram of exact sequences

$$
\CD
@.         0  @.   0   @. 0\\
@.   @VVV     @VVV        @VVV \\
0 @>>>   A_1 @>>>   \text {ker}(d)   @>>> \ker (\deg)=\oplus_{r=1}^s \Pic (H_r) \oplus (\oplus _{i=1}^t \Pic^0 (\widetilde E_i)) \\
@. @VVV        @VVV        @VVV\\
0@>>> A_2 @>>> \Pic (\Cal S_0) @>>> \Pic (S_0^{\text{nor}})= \oplus_{r=1}^s \Pic (H_r) \oplus (\oplus _{i=1}^t \Pic (\widetilde E_i))\\ 
@. @. @VdVV        @V\text {deg}VV\\
@. @.  M @= M    \\
\endCD
$$
where $A_1$ and $A_2$ are defined so that the above sequences are exact, and $A_2$ is finite as follows from the facts that
the sheaf $\rho _{*}(\Cal O_{S_0^{\text{nor}}} ^{\times})/\Cal O_{S_0}^{\times }$ is a skyscraper sheaf and the residue fields at closed points of $\Cal S_0$ are finite fields. The kernel 
$\text{ker}(\text {deg})=\oplus_{r=1}^s \Pic (H_r) \oplus (\oplus _{i=1}^t \Pic^0 (\widetilde E_i))$ of the right lower vertical map is finite: $\Pic^0 (\widetilde E_i)$ is finite since $\widetilde E_i$ is a proper and non-singular curve over a finite field, and for $1\le r\le s$ it holds $\Pic(H_r)$ is finite since $H_r$ is an affine and normal $1$-dimensional scheme of finite type over the finite field $F$. Indeed assume for simplicity that $H_r$ is geometrically connected over $F$.
Let $\ell/F$ be a finite extension such that $U_r\defeq H_r\times _{\Spec F}\Spec \ell$ admits a smooth and connected compactification $C_r$ with $(C_r\setminus U_r) (\ell )\neq \emptyset$. Let $U_r\to H_r$ be the canonical morphism and $\Pic(H_r)\to \Pic (U_r)$ the induced map of pull-back of line bundles. Then $\Ker [ \Pic(H_r)\to \Pic (U_r)]$ is finite (cf. [Guralnick-Jaffe-Raskind], Theorem 1.8). Further, the map $\Pic^0(C_r)\to  \Pic (U_r)$ obtained by restricting a degree $0$ line bundle on $C_r$ to $U_r$ is surjective [if $x\in (C_r \setminus U_r)(\ell)$ and $D\in \Pic (U_r)$ has degree $m$ then $D-mx\in \Pic^0(C_r)$ restricts to $D$ on $U_r$]
hence $\Pic (U_r)$ is finite since $\Pic ^0(C_r)$ is finite. From the above it follows that $\Pic(H_r)$ is finite.

This finishes the proof of Claim 1.
\qed
\enddemo

Consider the composite map
$$\psi_n:\Pic (\Cal S_n)\to \Pic (\Cal S_0) @>d>> M=\oplus_{i=1}^t(\oplus _{j=1}^{n_i}\Bbb Z).$$

\proclaim {Claim 2} $\text {ker}(\psi_n)$ is {\bf finite}
\endproclaim
\demo {Proof of Claim 2} First we prove that the kernel of the map
$\Pic (\Cal S_n)\to \Pic (\Cal S_{n-1})$ is finite for $n\ge 2$. Write $\Cal I_n$ for the sheaf of ideals of $\Cal O_{\Cal S}$ defining $\Cal S_n$. We have an exact sequence of sheaves on $\Cal S_n$:
$$1\to 1+(\Cal I_{n-1}/\Cal I_n)\to \Cal O_{S_n}^{\times}\to \Cal O_{S_{n-1}}^{\times}\to 1$$
which induces an exact sequence in cohomology
$$H^1(\Cal S_n,1+(\Cal I_{n-1}/\Cal I_n))\to \Pic(\Cal S_n)\to \Pic(\Cal S_{n-1})\to H^2(\Cal S_n,1+(\Cal I_{n-1}/\Cal I_n)).$$
Further the truncated exponential map $\alpha\mapsto 1+\alpha$ induces an isomorphism of sheaves $\Cal I_{n-1}/\Cal I_n\isom 1+(\Cal I_{n-1}/\Cal I_n)$ [$(\Cal I_{n-1}/\Cal I_n)^2=0$], hence $H^2(\Cal S_n,1+\Cal I_{n-1}/\Cal I_n)=0$ and the map $\Pic(\Cal S_n)\to \Pic(\Cal S_{n-1})$ is surjective. 
Moreover, $H^1(\Cal S_n,\Cal I_{n-1}/\Cal I_n)$ is finite. Indeed
$H^1(\Cal S_n,\Cal I_{n-1}/\Cal I_n)$ is a finitely generated $B_n$-module with finite support since the morphism $\lambda_n^{-1}(\Cal Z_n\setminus \{z_1,\ldots,z_t\})\to 
\Cal Z_n\setminus \{z_1,\ldots,z_t\}$ is affine and 
$R^1(\pi_n)_{*}(\Cal I_{n-1}/\Cal I_n)$ is the sheaf associated to the $B_n$-module $H^1(\Cal S_n,\Cal I_{n-1}/\Cal I_n)$, here $\lambda_n:\Cal S_n\to \Cal Z_n$ is the proper 
morphism induced by $\lambda$. 
This shows the kernel of the map
$\Pic (\Cal S_n)\to \Pic (\Cal S_{n-1})$ is finite for all $n\ge 2$. A similar argument shows that the kernel of the map $\Pic (\Cal S_1)\to \Pic (\Cal S_{0})$ is finite.
Hence, using Claim 1, $\text {ker}(\psi_n)$ is finite.

This finishes the proof of Claim 2.
\qed
\enddemo

In light of Claim 2, and in order to prove assertion {\bf (A)}, it suffices to prove that the cokernel of the composite map
$$ M \defeq \oplus_{i=1}^t(\oplus _{j=1}^{n_i}\Bbb Z) @>\beta>> \Pic(\Cal S)\to \Pic (\Cal S_n)  \to \Pic (\Cal S_0) @>d>> M=\oplus_{i=1}^t(\oplus _{j=1}^{n_i}\Bbb Z)$$
is finite. The latter follows from the nondegeneracy of the intersection pairing $(\oplus _{j=1}^{n_i}\Bbb Z)\times (\oplus _{j=1}^{n_i}\Bbb Z)\to \Bbb Z$ on each fibre
$E_i$ (cf. [Shafarevich], Lemma in page 69 and the discussion in page 71 after this Lemma), $1\le i\le t$.

This finishes the proof of assertion {\bf (A)}.
\qed
\enddemo

\demo{Proof of assertion {\bf (B)}} Let $\Cal J$ be an ample invertible $\Cal O_S$-ideal such that $\text{Supp} (\Cal O_S/\Cal J)=\Cal S_0$. The existence of such
$\Cal J$ follows from the facts that $H_r$ is affine (cf. Proof of Assertion A), $1\le r\le s$, the intersection pairing $(\oplus _{j=1}^{n_i}\Bbb Z)\times (\oplus _{j=1}^{n_i}\Bbb Z)\to \Bbb Z$ on each fibre $E_i$ is negative definite (cf. [Shafarevich], Lemma in page 69 and the discussion in page 71 after this Lemma), and the numerical criterion of ampleness on curves.
More precisely, $\forall \ 1\le i\le t$, one can find a divisor $D=\sum_{ j=1}^{n_i}m_{ij} D_{i,j}$ with $m_{i,j}<0$ and $D.D_{i,j} > 0$ for all $1\le j \le n_j$.

For $m\ge 1$, let $\Cal S'_m$ be the closed subscheme of $\Cal S$ defined by the sheaf of ideals $\Cal J^m$. To prove Assertion B it suffices to prove that there exists $m_0>0$ such that the map
$$\Pic (\Cal S'_{m+1})\to \Pic (\Cal S'_m)$$
is an isomorphism for any $m>m_0$. We have an exact sequence of shaves on $\Cal S'_{m+1}$:
$$1\to \Cal J^{m}/\Cal J^{m+1}\to \Cal O_{S'_{m+1}}^{\times}\to \Cal O_{S'_{m}}^{\times}\to 1$$
where the map $\Cal J^{m}/\Cal J^{m+1}\to \Cal O_{S'_{m+1}}^{\times}$ maps a  local section $\alpha$ to $1+\alpha$,
which induces an exact sequence in cohomology
$$H^1(\Cal S'_{m+1},\Cal J^{m}/\Cal J^{m+1})\to \Pic(\Cal S'_{m+1})\to \Pic(\Cal S'_{m})\to 0.$$
Now there exists $m_0>0$ such that $H^1(\Cal S'_{m+1},\Cal J^{m}/\Cal J^{m+1})=0$ if $m\ge m_0$ by [EGA III], premi\`ere partie, Proposition 2.2.1.

This finishes the proof of assertion {\bf (B)}.
\qed
\enddemo

This finishes the proof of Proposition 4.1.
\qed

\subhead
\S 5. Compactification of formal germs of $p$-adic curves
\endsubhead
In this section we use the following notations: $K$ is a complete discrete valuation field with valuation ring $R$, uniformising parameter $\pi$,  
and with perfect residue field $\ell\defeq R/\pi R$. Further, $A$ is a {\bf two dimensional normal complete local ring}
containing $R$ with maximal ideal $\frak m_A$ containing $\pi$ and residue field $\ell=A/{\frak m}_A$. 
We assume that $X\defeq \Spec (A\otimes_RK)$ is geometrically connected. Given a finite extension $L/K$ we write $\Cal O_L$ for the valuation ring of $L$, 
$A_L\defeq A\otimes _{\Cal O_L}L$, $A_{\Cal O_L}\defeq A\otimes _R{\Cal O_L}$, and $A_{\Cal O_L}^{\nor}$  the normalisation of $A_{\Cal O_L}$ in its total ring of fractions.

\proclaim {Proposition 5.1 (Compactification of formal germs of $p$-adic curves)} 
We use the above notations. There exists a finite extension $L/K$, a flat, proper, connected, and normal $\Cal O_L$-relative curve
$\Cal Y\to \Spec \Cal O_L$, a closed point $y\in \Cal Y$, and an isomorphism $\hat \Cal O_{\Cal Y,y}\isom A_{\Cal O_L}^{\nor}$ where $\hat \Cal O_{\Cal Y,y}$ is the completion of the local ring $\Cal O_{\Cal Y,y}$ of $\Cal Y$ at $y$.
\endproclaim

\demo{Proof} By the main result in [Epp], Introduction,  there exists a finite extension $L/K$ with uniformising parameter $\pi_L$ such that 
$A_{\Cal O_L}^{\nor}/\pi_LA_{\Cal O_L}^{\nor}$ 
is reduced. Note that $A_{\Cal O_L}^{\nor}$ is a normal two dimensional complete local ring with perfect residue field (cf. [Bourbaki], Chapter IX, $\S4$, Lemma 1, and our assumption that $X$ is geometrically connected).
Without loss of generality we will assume that $A/\pi A$ is reduced. We show there exists a proper, flat, connected, and normal relative $R$-curve
$\Cal Y\to \Spec R$, a closed point $y\in \Cal Y$, and an isomorphism $\hat \Cal O_{\Cal Y,y}\isom  A$. 

First, $A/\pi A$ is a (reduced) one dimensional complete local ring with residue field $\ell$, hence is isomorphic to a quotient
$\ell[[x_1,\ldots,x_t]]/{\frak a}$ of a formal power series ring $\ell[[x_1,\ldots,x_t]]$
over $\ell$ (cf. [Bourbaki], Chapitre IX, $\S3$).  It then follows from [Artin], Theorem 3.8, and basic facts on the theory of algebraic curves, that there exists a proper and reduced 
connected (but not necessarily irreducible) $\ell$-curve $Z$, a closed point $y\in Z$, and an isomorphism $\hat \Cal O_{Z,y}\isom A/\pi A$ where
$\hat \Cal O_{Z,y}$ is the completion of the local ring $\Cal O_{Z,y}$ of $Z$ at $y$. Moreover, $Z$ is non-singular outside $y$.
There exists a rational function $f$ on $Z$ which defines a finite generically separable morphism $f:Z\to \Bbb P^1_{\ell}$
such that $y=f^{-1}(\infty)$ (cf. [Harbater-Stevenson], proof of Theorem 3). Thus, by considering the completion of the morphism $f$ above $\infty$, we obtain a finite generically separable 
morphism $\bar g:\Spec (A/\pi A)\to \Spec (\ell[[t]])$ where $t$ is a local parameter at $\infty$.
This morphism lifts to a finite morphism $g:\Spf A\to \Spf (R[[T]])$ of formal schemes (cf. loc. cit., Lemma 2).
Let $\widetilde Z\to Z$ be the morphism of normalisation and $\{x_1,\ldots,x_m\}\subset \widetilde Z$ the pre-image of $y$. 
There is a one-to-one correspondence between the set $\{\frak p_1,\ldots,\frak p_m\}\subset \Spec A$ of prime ideals of height $1$ containing $\pi$ and the set $\{x_1,\ldots,x_m\}$, 
$\frak p_i$ corresponds to $x_i$, $1\le i\le m$. The composite morphism $\widetilde Z\to Z\to \Bbb P^1_{\ell}$ induces, by completion above $\infty$, finite separable morphisms 
$\bar g_i:\Spec \Fr (\hat \Cal O_{\widetilde Z,x_i})\to \Spec \ell((t))$ where $\Fr (\hat \Cal O_{\widetilde Z,x_i})$ is the fraction field of the completion $\hat \Cal O_{\widetilde Z,x_i}$
of the local ring $\Cal O_{\widetilde Z,x_i}$ of $\widetilde Z$ at $x_i$, $1\le i\le m$ (with the above notations $t=T \mod \pi$).

Consider the formal closed unit disc $D=\Spf R<\frac {1}{T}>$ with parameter $\frac {1}{T}$ and its special fibre $D_{\ell}=\Spec \ell [\frac {1}{t}]$ ($D_{\ell} \isom \Bbb A^1_{\ell}$).
By a result of Gabber and Katz (cf. [Katz], Main Theorem 1.4.1) there exists, for $1\le i\le m$, a finite cover $\bar h_i:C_i\to D_{\ell}$ with $C_i$ connected, which only (tamely) ramifies above the point $\frac {1}{t}=0$
and such that the completion of $\bar h_i$ above $t=0$ is generically isomorphic to the cover $\bar g_i:\Spec \Fr (\hat \Cal O_{\widetilde Z,x_i})\to \Spec \ell ((t))$.
Using formal patching techniques (cf. [Sa\"\i di4], 1.2) one can lift the covers $\bar h_i$ to finite covers $h_i:Y_i\to D$ which only ramify 
above the point $\frac {1}{T}=0$, $1\le i\le m$ [outside $\frac {1}{T}=0$
the existence of such a lifting follows from the theorems of lifting of \'etale covers (cf. [Grothendieck], Expos\'e I, Corollaire 8.4). In a formal neighbourhood of $\frac {1}{T}=0$ such a lifting is possible under the tameness condition: \'etale locally near $\frac {1}{t}$ the cover $\bar h_i$ is defined by an equation $y^s=\frac{1}{t^e}$, where $s\ge 1$ is an integer prime to the 
characteristic of $\ell$, and one lifts to the cover defined by $Y^s=\frac{1}{T^e}$].
For $1\le i\le m$, let $\hat A_{\frak p_i}$ be the completion of the localisation $A_{\frak p_i}$ of $A$ at ${\frak p}_i$. Thus, 
$\hat A_{\frak p_i}$ is a complete discrete valuation ring with uniformising parameter $\pi$ (recall $A/\pi A$ is reduced) 
and residue field $\Fr (\hat \Cal O_{\widetilde Z,x_i})$. Let $B$ be the completion of the localisation of $R[[T]]$ at $\pi$. Thus, $B$ is a complete discrete valuation ring with residue field 
$\ell((t))$. The finite cover $g:\Spf A\to \Spf (R[[T]])$ induces, by pull-back to $\Spf B$, finite covers $g_i:\Spf \hat A_{\frak p_i}\to \Spf B$ which (by construction) lift the covers 
$\bar g_i:\Spec \Fr (\hat \Cal O_{\widetilde Z,x_i})\to \Spec \ell((t))$, $1\le i\le m$. Further, the cover $h_i:Y_i\to D$ induces, by pull-back to $\Spf B$, a finite cover
$\tilde h_i:\Spf B_i\to \Spf B$ which by construction lifts the cover $\bar g_i:\Spec \Fr (\hat \Cal O_{\widetilde Z,x_i})\to \Spec \ell((t))$. Thus, the covers $\tilde h_i:\Spf B_i\to \Spf B$
and $g_i:\Spf \hat A_{\frak p_i}\to \Spf B$ are isomorphic since $\bar g_i$ is generically separable. 
Using formal patching techniques (cf. loc. cit.) one can patch the covers $g:\Spf A\to \Spf (R[[T]])$ and 
$h_i:Y_i\to D$, $1\le i\le m$, to construct a finite cover $\Cal Y\to \Bbb P^1_{R}$ in the category of formal schemes with $\Cal Y$ normal, connected, proper, and flat over $\Spf R$.
The special fibre $\Cal Y_{\ell}\defeq \Cal Y\times _{\Spec R}\Spec \ell$ of $\Cal Y$ consists of $m$ irreducible components which intersect at the point $y$ and is (by construction) non-singular outside $y$. The formal curve $\Cal Y$ is algebraic by formal GAGA and (by construction) $\hat \Cal O_{\Cal Y,y}\isom A$ as required.
\qed
\enddemo

\definition {Remark 5.2} Proposition 5.1 asserts the existence, after possibly a finite extension of $K$, of a proper $R$-curve $\Cal Y$ and a closed point $y\in \Cal Y^{\cl}$ such that  
$\hat \Cal O_{\Cal Y,y}\isom A$. The special fibre $\Cal Y_{\ell}\defeq \Cal Y\times _{\Spec R}\Spec \ell$ of $\Cal Y$ consists of $m_y\defeq m$ (cf. the proof of Proposition 5.1 for the definition of $m$) irreducible components $\{C_1,\ldots,C_m\}$ which intersect at $y$, $\Cal Y_{\ell}$ is non-singular outside $y$, and the normalisation morphism $C_i^{\nor}\to C_i$ is a homeomorphism, $1\le i\le m$. 
In fact one can, assuming the existence of a compactification of $\Spec A$ as in Proposition 5.1, construct such a compactification $\Cal Y$ of $\Spec A$ with the additional property that $C_i^{\nor}\isom \Bbb P^1_{\ell}$, $\forall 1\le i\le m$ (cf. [Sa\"\i di4], Remark 3.1).
\enddefinition

\proclaim {Proposition 5.3} We use the above notations. There exists a finite extension $L/K$ and a finite morphism $\Spec B\to \Spec A_{\Cal O_L}^{\nor}$ with $B$ local, normal,
{\bf hyperbolic} (cf. Notations), and the morphism $\Spec B_L\to \Spec A_L$ is geometric and \'etale.
\endproclaim

\demo{Proof} This follows easily from Proposition 5.1, Remark 5.2,  and Theorem 3 in [Sa\"\i di4].
\qed
\enddemo

$$\text{References.}$$

\noindent
[Artin] Artin, M., Algebraic approximation of structures over complete local rings,
Publications math\'ematiques de l'I.H.\'E.S., tome 36 (1969), p. 23-58.

\noindent
[Bourbaki] Bourbaki, N., Alg\`ebre Commutative, Chapitre 9, Masson, 1983.

\noindent
[Cossart-Jannsen-Saito], Cossart, V., Jannsen, U., Saito, S., Canonical embedded and non-embedded resolution of singularities for excellent two-dimensional schemes,
arXiv:0905.2191.

\noindent
[EGA III] Grothendieck, A., Dieudonn\'e, J., \'Etude cohomologique des faisceaux coh\'erents, Pub. Math. IHES 11 (1961), 17 (1963). 

\noindent
[Epp] Epp, P., H., Eliminating wild ramification, Inventiones math. 19 (1973), 235- 249.

\noindent
[Esnault-Wittenberg] Esnault, H., Wittenberg, O., On abelian birational sections, Journal of the American Mathematical Society, 
Volume 23, Number 3 (2010), 713-724.

\noindent
[Fresnel-Matignon] Fresnel, J., Matignon, M., Sur les espaces analytiques quasi-compacts de dimension 1 sur un corps valu\'e, complet, ultram\'etrique, Ann. Mat. Pura. Appl., (IV), Vol. CXLV (1986), 159-210.

\noindent
[Gabber-Liu-Lorenzini] Gabber, O., Liu, Q., Lorenzini, D., The index of an algebraic variety, Invent. Math. 192 (2013), no. 3, 567-626.

\noindent
[Grothendieck] Grothendieck, A., Rev\^etements \'etales et groupe fondamental, Lecture 
Notes in Math. 224, Springer, Heidelberg, 1971.

\noindent
[Guralnick-Jaffe-Raskind] Guralnick, R., Jaffe, D.B., Raskind, W., On the Picard Group: Torsion and the Kernel Induced by a Faithfully Flat Map, Journal of Algebra 183 (1996), 420-455.

\noindent
[Hansen] Hansen, D., Vanishing and comparison theorems in rigid analytic geometry, Compositio Mathematica,
Vol. 156, Issue 2 , (2020) , pp. 299-324.

\noindent
[Harbater-Stevenson] Harbater, D., Stevenson, K., Patching and thickening problems, 
Journal of Algebra 212 (1999), 272-304.

\noindent
[Hoshi] Hoshi, Y., Existence of nongeometric pro-$p$ Galois sections of hyperbolic curves, Publ. Res. Inst. Math. Sci. 46 (2010), no. 4, 829-848.

\noindent
[Katz] Katz, N., Local-to-global extensions of representations of fundamental groups, Ann.
Inst. Fourier (Grenoble) 36 (1986), no. 4, 69-106.

\noindent
[Lipman] Lipman, J., Desingularization of two dimensional schemes, Ann. of Math. 107 (1978), 151-207.

\noindent
[Liu] Liu, Q., Algebraic geometry and arithmetic curves, Oxford graduate texts in mathematics 6. Oxford University Press, 2002.


\noindent
[Mochizuki] Mochizuki, S., Absolute anabelian cuspidalizations of proper hyperbolic curves,  J. Math. Kyoto
Univ.  47  (2007),  no. 3, 451--539.

\noindent
[Neukirch-Schmidt-Wingberg] Neukirch, J., Schmidt, A., Wingberg, K., Cohomology of number fields, first edition,
Springer, Grundlehren der mathematischen Wissenschaften Bd. 323, 2000.

\noindent
[Sa\"\i di1] Sa\"\i di, M., The cuspidalisation of sections of arithmetic fundamental groups, Advances in Mathematics 230
(2012), 1931-1954.

\noindent
[Sa\"\i di2]  Sa\"\i di, M., The cuspidalisation of sections of arithmetic fundamental groups II, Advances in Mathematics 345
(2019), https://doi.org/10.1016/j.aim.2019.106737.

\noindent
[Sa\"\i di3] Sa\"\i di, M., \'Etale fundamental groups of affinoid $p$-adic curves, Journal of algebraic geometry, 27 (2018), 727--749.

\noindent
[Sa\"\i di4] Sa\"\i di, M., On \'etale fundamental groups of formal germs of $p$-adic curves, Tohoku Math journal, volume 72, no. 1 (2020), pages 63-76.

\noindent
[Sa\"\i di5]  Sa\"\i di, M., Arithmetic of $p$-adic curves and sections of geometrically abelian fundamental groups. Math. Z. 297, no. 3-4 (2021), 1191-1203.

\noindent
[Sa\"\i di6] Sa\"\i di, M., On the existence of non-geometric sections of arithmetic fundamental groups, 
Math. Z. 277, no. 1-2 (2014), 361-372.


\noindent
[Shafarevich] Shafarevich, I.R., Lectures on minimal models and birational transformations of two dimensional schemes. Tata Inst. Fund. Research 37, Bombay, 1966.

\noindent
[Serre] Serre J.-P., Construction de rev\^etements de la droite affine en caract\'eristique $p>0$, C. R. Acad. Sci. Paris 311 (1990), 341-346.

\noindent
[Tamagawa] Tamagawa, A., The Grothendieck conjecture for affine curves,  Compositio Math.  109  (1997),  no. 2, 135--194.

\bigskip

\noindent
Mohamed Sa\"\i di

\noindent
College of Engineering, Mathematics, and Physical Sciences

\noindent
University of Exeter

\noindent
Harrison Building

\noindent
North Park Road

\noindent
EXETER EX4 4QF

\noindent
United Kingdom

\noindent
M.Saidi\@exeter.ac.uk

\end
\enddocument